\documentclass[a4paper, 12pt]{article}
\usepackage{enumerate, theorem}
\usepackage{amsmath, amsfonts, amssymb}
\usepackage[height=22.5cm, width=15cm]{geometry}
\usepackage[all]{xy}
\usepackage{mathrsfs, yfonts}
\newcommand{\A}{\tilde{\mathcal{A}}}

\newcommand{\parag}[1]{\paragraph{\sc{#1.}}}

\newtheorem{thm}{Th\'eor\`eme}[subsection]
\newtheorem{defn}[thm]{D\'efinition}
\newtheorem{cor}[thm]{Corollaire}
\newtheorem{prop}[thm]{Proposition}
\newtheorem{lemma}[thm]{Lemme}

\begin{document}

\title{Two finiteness theorem for (a,b)-modules.}

\date{16/01/08} 
 
\author{Daniel Barlet\footnote{Barlet Daniel, Institut Elie Cartan UMR 7502  \newline
Nancy-Universit\'e, CNRS, INRIA  et  Institut Universitaire de France, \newline
BP 239 - F - 54506 Vandoeuvre-l\`es-Nancy Cedex.France.\newline
e-mail : barlet@iecn.u-nancy.fr}.}
 
 \maketitle
 
 \markright{Filtered Gauss-Manin connexion.}
 
 \section*{Summary.}
 
 We prove the following two results
 \begin{enumerate}
 \item  For a proper holomorphic function \ $ f : X \to D$ \ of a complex manifold \ $X$ \ on a disc such that \ $\{ df = 0 \} \subset f^{-1}(0)$, we construct, in a functorial  way, for each integer \ $p$,  a geometric (a,b)-module \ $E^p$ \ associated to the (filtered) Gauss-Manin connexion of \ $f$.\\
  This first theorem is an existence/finiteness result which shows that geometric (a,b)-modules may be used in  global situations.
  \item For any regular (a,b)-module \ $E$ \ we give an integer \ $N(E)$, explicitely given from simple invariants of \ $E$, such that the isomorphism class of \ $E\big/b^{N(E)}.E$ \ determines the isomorphism class of \ $E$.\\
  This second result allows to cut asymptotic expansions (in powers of \ $b$) \ of elements of \ $E$ \ without loosing any information.
\end{enumerate}

\bigskip

AMS Classification :  32-S-05,  32S 20, 32-S-25, 32-S-40.

\bigskip

Key words : (a,b)-module or Brieskorn modules, Gauss-Manin connexion, vanishing cycles.

\newpage

 \tableofcontents
 
 \newpage
 
 \section{Introduction.}
 
 The following situation is frequently met : we consider a vector space \ $E$ \ of multivalued holomorphic functions (possibly with values in a complex finite dimensional vector space \ $V$) with finite determination on a punctured disc around \ $0$ in \ $\mathbb{C}$, stable by multiplication by the variable \ $z$ \ and by "primitive". These functions are determined by there formal asymptotic expansions at \ $0$ \ of the type
  $$ \sum_{(\alpha,j)\in A\times[0,n]}  c_{\alpha,j}(z).z^{\alpha}.\frac{(Log z)^j}{j!}  $$
  where \ $n$ \ is a fixed integer, where \ $A$ \ is a finite set of complex numbers whose real parts are, for instance, in the interval \ $ ]-1,0]$, and where  the  \ $c_{\alpha,j}$ \ are in \ $\mathbb{C}[[z]] \otimes V$.\\
  Define  \ $e(\alpha,j) = z^{\alpha}.(Log z)^j/j! $ \ and \ $E(\alpha,n) : = \oplus_{j=0}^n \ (\mathbb{C}[[z]]\otimes V).e(\alpha,j)$. Then for each \ $\alpha, \,E(\alpha,n)$ \ is a free \ $\mathbb{C}[[z]]-$module of rank \ $(n+1).\dim V$ \ which is stable by \ $b : = \int_0^z $. To be precise, \ $b$ \ is defined by induction on \ $j \geq 0$ \ by the "obvious" formulas :
   $$b[e(\alpha,0)] = \frac{e(\alpha+1,0)}{\alpha+1} \quad {\rm and \  for} \quad j \geq 1$$
  $$b[e(\alpha,j)] = \frac{e(\alpha+1,j)}{\alpha+1} -\frac{1}{\alpha+1}.b[e(\alpha,j-1)].$$
  So we have an inclusion \ $ E \subset \oplus_{\alpha \in A} E(\alpha)$ \ which is compatible with \ $a : = \times z$ \ and by  \ $b$. Note that each \ $E(\alpha)$ \  is also a free finite rank  \ $\mathbb{C}[[b]]-$module which is stable by \ $a$.\\
  Let us assume now that \ $E$ \ is also a  \ $\mathbb{C}[[b]]-$module which is stable by \ $a$. Then \ $E$ \ has to be free and of finite rank over \ $\mathbb{C}[[b]]$. Our aim in this situation is to understand and to describe the relations between the coefficients \ $c_{\alpha,j}$ \ of the asymptotic expansions of elements in \ $E$. This leads to construct "invariants" associated to the given \ $E$.\\
  As in the "geometric" situations we consider the complex numbers \ $\alpha \in A$ \ are rationnal numbers, a change of variable of type \ $t : = z^{1/N}$, with \ $N \in \mathbb{N}^*$, allows to reduce the situation to the case where all \ $\alpha$ \ are \ $0$. Then one can use the "nilpotent" operator \ $\frac{d}{dt}$ \ and the nilpotent part of the monodromy  \ $e_{\alpha,j} \to e_{\alpha,j-1}$ \ to construct some filtrations in order, in the good cases, to build a Mixte Hodge structure on a finite dimensional vector space associated to \ $E$. For instance, this is the case in A.N. Varchenkho's description[V. 80] of the Mixte Hodge structure built by J. Steenbrink [St. 76] on the cohomology of the Milnor fiber of an holomorphic function with an isolated singularity at the origin of \ $\mathbb{C}^{n+1}$.\\
  
  \smallskip
  
  But it is clear that we loose some information in this procedure. The point of view which is to consider \ $E$ \ itself  as a left module on the (non commutative) algebra
   $$ \A : = \{ \sum_{\nu \geq 0} \ P_{\nu}(a).b^{\nu} \}$$
   where the \ $P_{\nu}$ \ are polynomials with complex coefficients is richer. This is evidenced by M. Saito result  [Sa. 91].
   
   \smallskip
   
   The aim of the first part of this article is to build in a natural way, for any proper holomorphic function \ $f : X \to D$ \ of a complex manifold \ $X$, assumed to be smooth outside its \ $0-$fiber \ $X_0 : = f^{-1}(0)$, a regular (geometric) (a,b)-module for each degree \ $p \geq 0$, which represent a filtered version of the Gauss-Manin connexion of \ $f$ \ at the origin.\\
   This result is in fact a finiteness theorem which is a first step to refine the limite Mixte Hodge structure in this situation. It is interesting to remark that no K{\"a}hler assumption is used in this construction of these geometric (a,b)-modules.\\
   This obviuosly shows that  (a,b)-modules are basic objects and that they are important not only in the study of  local singularities of holomorphic functions but more generally  in complex geometry . So it is interesting to have some tools in order  to compute them.\\
   This is precisely the aim of the second part of this paper. We prove a finiteness result which gives, for a regular (a,b)-module \ $E$, an integer \ $N(E)$, bounded by simple numerical invariants of \ $E$, such that you may cut the asymptotic expansions (in powers of \ $b$) of elements of \ $E$ \ without any lost of information on the structure of the (a,b)-module \ $E$. It is well known that the formal asymptotic expansions for solutions of a regular differential system always converge, and also that such an integer exists for any meromorphic connexion in one variable (see [M.91] proposition 1.12). But it is important to have an effective bound for such an integer  easely computable from simple invariants of the (a,b)-module structure of \ $E$.

  \section{The existence theorem.}
 
 \subsection{Preliminaries.}
 
 Here we shall complete and precise the results of the section 2 of [B.07]. The situation we shall consider is the following : let \ $X$ \ be a connected  complex manifold of dimension \ $n +1$ \ and \ $f : X \to \mathbb{C}$ \ a non constant holomorphic function such that \ $\{ x \in X / \ df = 0 \} \subset f^{-1}(0)$. We introduce the following complexes of sheaves supported by \ $X_0 : = f^{-1}(0)$
 \begin{enumerate}
 \item  The formal completion ''in \ $f$'' \ $(\hat{\Omega}^{\bullet}, d^{\bullet})$ \ of the usual holomorphic de Rham complex of \ $X$.
 \item The sub-complexes \ $(\hat{K}^{\bullet}, d^{\bullet})$ \ and \ $(\hat{I}^{\bullet}, d^{\bullet})$ \ of \  $(\hat{\Omega}^{\bullet}, d^{\bullet})$ \  where the subsheaves \ $\hat{K}^p$ \ and \ $\hat{I}^{p+1}$ \ are defined for each \ $p \in \mathbb{N}$ \  respectively as the kernel and the image of the map
 $$  \wedge df : \hat{\Omega}^p \to \hat{\Omega}^{p+1} $$
 given par exterior multiplication by \ $df$.  We have the exact sequence
 \begin{equation*} 0 \to (\hat{K}^{\bullet}, d^{\bullet}) \to (\hat{\Omega}^{\bullet},  d^{\bullet}) \to (\hat{I}^{\bullet}, d^{\bullet})[+1]  \to 0. \tag{1}
 \end{equation*}
 Note that \ $\hat{K}^0$ \ and \ $\hat{I}^0$ \ are zero by definition.
 \item The natural inclusions \ $\hat{I}^p \subset \hat{K}^p$ \ for all \ $p \geq 0$ \ are compatible with the differential \ $d$. This leads to an exact sequence of complexes
 \begin{equation*}
 0 \to (\hat{I}^{\bullet}, d^{\bullet}) \to (\hat{K}^{\bullet}, d^{\bullet}) \to ([\hat{K}/\hat{I}]^{\bullet}, d^{\bullet}) \to 0 .\tag{2}
 \end{equation*}
 \item We have a natural inclusion \ $f^*(\hat{\Omega}_{\mathbb{C}}^1) \subset \hat{K}^1\cap Ker\, d$, and this gives a sub-complex (with zero differential) of \ $(\hat{K}^{\bullet}, d^{\bullet})$. As in [B.07], we shall consider also the complex \ $(\tilde{K}^{\bullet}, d^{\bullet})$ \ quotient. So we have the exact sequence
 \begin{equation*}
  0 \to f^*(\hat{\Omega}_{\mathbb{C}}^1) \to (\hat{K}^{\bullet}, d^{\bullet}) \to (\tilde{K}^{\bullet}, d^{\bullet}) \to 0 . \tag{3}
  \end{equation*}
 We do not make the assumption here that \ $f = 0 $ \ is a reduced equation of \ $X_0$, and we do not assume that \ $n \geq 2$, so the cohomology sheaf in degree 1  of the complex \ $(\hat{K}^{\bullet}, d^{\bullet})$, which is equal to \ $\hat{K}^1 \cap Ker\, d$ \ does not co{\"i}ncide, in general, with \ $f^*(\hat{\Omega}_{\mathbb{C}}^1)$. So the complex \ $ (\tilde{K}^{\bullet}, d^{\bullet})$ \ may have a non zero cohomology sheaf in degree 1.
  \end{enumerate} 
  Recall now that we have on the cohomology sheaves of the following complexes \\
   $(\hat{K}^{\bullet}, d^{\bullet}), (\hat{I}^{\bullet}, d^{\bullet}), ([\hat{K}/\hat{I}]^{\bullet}, d^{\bullet}) $ \ and \ $f^*(\hat{\Omega}_{\mathbb{C}}^1), (\tilde{K}^{\bullet}, d^{\bullet})$ \ natural operations \ $a$ \ and \ $b$ \ with the relation \ $a.b - b.a = b^2$. They are defined in a na{\"i}ve way by 
  $$  a : = \times f \quad {\rm and} \quad  b : = \wedge df \circ d^{-1} .$$
  The definition of \ $a$ \ makes sens obviously. Let me precise the definition of \ $b$ \ first in  the case of \ $\mathcal{H}^p(\hat{K}^{\bullet}, d^{\bullet})$ \ with \ $p \geq 2$ \ : if \ $x \in \hat{K}^p \cap Ker\, d$ \ write \ $x = d\xi$ \ with \ $\xi \in \hat{\Omega}^{p-1}$ \ and let \ $b[x] : = [df\wedge \xi]$. The reader will check easily that this makes sens.\\
  For \ $p = 1$ \ we shall choose \ $\xi \in \hat{\Omega}^0$ \ such that \ $\xi = 0$ \ on the smooth part of  \ $X_0$ \ (set theoretically). This is possible because the condition \ $df \wedge d\xi = 0 $ \ allows such a choice : near a smooth point of \ $X_0$ \ we can choose coordinnates such \ $ f = x_0^k$ \ and the condition on \ $\xi$ \ means independance of \ $x_1, \cdots, x_n$. Then \ $\xi$ \ has to be (set theoretically) locally constant on \ $X_0$ \ which is locally connected. So we may kill the value of such a \ $\xi$ \ along \ $X_0$.\\
  The case of the complex \ $(\hat{I}^{\bullet}, d^{\bullet})$ \ will be reduced to the previous one using the next  lemma.
  
  \begin{lemma}\label{tilde b}
  For each \ $p \geq 0$ \ there is a natural injective map
  $$\tilde{b} :  \mathcal{H}^p(\hat{K}^{\bullet}, d^{\bullet}) \to \mathcal{H}^p(\hat{I}^{\bullet}, d^{\bullet})$$
  which satisfies the relation \ $a.\tilde{b} = \tilde{b}.(a + b) $. For \ $p \not= 1$ \ this map is bijective.
  \end{lemma}
  
  \parag{Proof} Let \ $x \in \hat{K}^p \cap Ker\, d $ \ and write \ $x = d\xi $ \ where \ $\xi \in \hat{\Omega}^{p-1}$ \ (with \ $\xi = 0$ \ on \ $X_0$ \ if \ $p = 1$), and set \ $\tilde{b}([x]) : = [df\wedge \xi] \in \mathcal{H}^p(\hat{I}^{\bullet}, d^{\bullet})$. This is independant on the choice of \ $\xi$ \ because, for \ $p \geq 2$, adding \ $d\eta$ \ to \ $\xi$ \ does not modify the result as \ $[df\wedge d\eta] = 0 $. For \ $p =1$ \ remark that our choice of \ $\xi$ \ is unique.\\
   This is also independant of the the choice of \ $x $ \ in \ $ [x] \in \mathcal{H}^p(\hat{K}^{\bullet}, d^{\bullet})$ \  because adding \ $\theta \in \hat{K}^{p-1}$ \ to \ $\xi$ \ does not change \ $df \wedge \xi$.\\
   Assume \ $\tilde{b}([x]) = 0 $ \ in \ $ \mathcal{H}^p(\hat{I}^{\bullet}, d^{\bullet})$; this means that we may find \ $\alpha \in \hat{\Omega}^{p-2}$ \ such \ $df \wedge \xi = df \wedge d\alpha$. But then, \ $\xi - d\alpha $ \ lies in \ $\hat{K}^{p-1}$ \ and \ $x = d(\xi - d\alpha ) $ \ shows that \ $[x] = 0$. So \ $\tilde{b}$ \ is injective.\\
  Assume now \ $p \geq 2$.  If \ $df\wedge \eta $ \ is in \ $\hat{I}^p \cap Ker\, d$, then \ $df \wedge d\eta = 0 $ \ and \ $y : = d\eta $ \ lies in \ $\hat{K}^p \cap Ker\, d$ \ and defines a class \ $[y] \in  \mathcal{H}^p(\hat{K}^{\bullet}, d^{\bullet}) $ \ whose image by \ $\tilde{b}$ \ is \ $[df\wedge \eta] $. This shows the surjectivity of \ $\tilde{b}$ \ for \ $p \geq 2$.\\
   For \ $p=1$ \ the map \ $\tilde{b}$ \ is not surjective (see the remark below).\\
  To finish the proof let us  to compute \ $\tilde{b}(a[x] + b[x])$. Writing again \ $x = d\xi$, we get
   $$ a[x] + b[x] =[ f.d\xi + df \wedge \xi] = [d(f.\xi)] $$
   and so
   $$ \tilde{b}( a[x] + b[x] ) = [ df \wedge f.\xi ] = a.\tilde{b}([x]) $$
   which concludes the proof. $\hfill \blacksquare$
   
   \bigskip
   
   Denote by \ $i :  (\hat{I}^{\bullet}, d^{\bullet}) \to (\hat{K}^{\bullet}, d^{\bullet})$ \ the natural inclusion and define the action of \ $b$ \ on \ $\mathcal{H}^p(\hat{I}^{\bullet}, d^{\bullet})$ \ by \ $b : = \tilde{b}\circ \mathcal{H}^p(i) $. As \ $i$ \ is \ $a-$linear, we deduce the relation \ $a.b - b.a = b^2$ \ on \ $\mathcal{H}^p(\hat{I}^{\bullet}, d^{\bullet})$ \ from the relation of the previous lemma. \\
   
   The action of \ $a$ \ on the complex \ $ ([\hat{K}/\hat{I}]^{\bullet}, d^{\bullet}) $ \ is obvious and the action of \ $b$ \ is zero.\\
   
   The action of \ $a$ \ and \ $b$ \ on \ $f^*(\hat{\Omega}_{\mathbb{C}}^1) \simeq E_1\otimes \mathbb{C}_{X_0}$ \ are the obvious one, where \ $E_1$ \ is the rank 1 (a,b)-module with generator \ $e_1$ \ satisfying \ $a.e_1 = b.e_1$ \ (or, if you prefer, \ $E_1 : = \mathbb{C}[[z]]$ \ with \ $a : = \times z,\quad b : = \int_0^z $ \ and \ $e_1 : = 1$). \\
   Remark that the natural inclusion \ $f^*(\hat{\Omega}^1_{\mathbb{C}}) \hookrightarrow (\hat{K}^{\bullet}, d^{\bullet})$ \ is compatible with the actions of \ $a$ \ and \ $b$. The actions of \ $a$ \ and \ $b$ \ on \ $\mathcal{H}^1(\tilde{K}^{\bullet}, d^{\bullet}) $ \ are simply  induced by the corresponding actions on \ $\mathcal{H}^1(\hat{K}^{\bullet}, d^{\bullet})$.
   
   \parag{Remark} The exact sequence of complexes (1) induces  for any \ $p \geq 2$ \  a bijection
   $$ \partial^p : \mathcal{H}^p(\hat{I}^{\bullet}, d^{\bullet}) \to \mathcal{H}^p(\hat{K}^{\bullet}, d^{\bullet})$$
   and a short exact sequence 
   \begin{equation*}
    0 \to \mathbb{C}_{X_0} \to \mathcal{H}^1(\hat{I}^{\bullet}, d^{\bullet}) \overset{\partial^1}{\to} \mathcal{H}^1(\hat{K}^{\bullet}, d^{\bullet}) \to 0 \tag{@}
    \end{equation*}
   because of the de Rham lemma. Let us check  that for \ $p \geq 2$ \ we have \ $\partial^p = (\tilde{b})^{-1}$ \ and that for \ $p =1$ \ we have \ $\partial^1\circ \tilde{b} = Id$. If \ $x = d\xi \in \hat{K}^p \cap Ker\, d$ \ then \ $\tilde{b}([x]) = [df\wedge \xi]$ \ and \ $\partial^p[df\wedge\xi] = [d\xi]$. So \ $ \partial^p\circ\tilde{b} = Id \quad \forall p \geq 0$. For \ $p \geq 2$ \ and \ $df\wedge\alpha \in \hat{I}^p \cap Ker\, d$ \ we have \ $\partial^p[df\wedge\alpha] = [d\alpha]$ \ and \ $\tilde{b}[d\alpha] = [df\wedge\alpha]$, so \ $\tilde{b}\circ\partial^p = Id$. For \ $p = 1$ \ we have \ $\tilde{b}[d\alpha] = [df\wedge(\alpha - \alpha_0)]$ \ where \ $\alpha_0 \in \mathbb{C}$ \ is such that \ $\alpha_{\vert X_0} = \alpha_0$.  This shows that in degree 1 \ $\tilde{b}$ \ gives a canonical splitting of the exact sequence \ $(@)$.
   
  \subsection{$\A-$structures.}
   
   Let us consider now the \ $\mathbb{C}-$algebra
   $$ \A : = \{ \sum_{\nu \geq  0} \quad P_{\nu}(a).b^{\nu} \}$$
   where \ $P_{\nu} \in \mathbb{C}[z]$, and the commutation relation \ $a.b - b.a = b^2$, assuming that left and right multiplications by \ $a$ \ are continuous for the \ $b-$adic topology of \ $\A$.
   
   Define the following complexes of sheaves of left  \ $\A-$modules on \ $X$ :
   \begin{align*}
   & ({\Omega'}^{\bullet}[[b]], D^{\bullet})  \quad {\rm and} \quad   ({\Omega''}^{\bullet}[[b]], D^{\bullet}) \quad {\rm where} \tag{4}\\
    & {\Omega'}^{\bullet}[[b]] : = \sum_{j=0}^{+\infty} b^j.\omega_j \quad {\rm with} \quad \omega_0 \in \hat{K}^p  \\
     & {\Omega''}^{\bullet}[[b]] : = \sum_{j=0}^{+\infty} b^j.\omega_j \quad {\rm with} \quad \omega_0 \in \hat{I}^p \\
     & D(\sum_{j=0}^{+\infty} b^j.\omega_j) = \sum_{j=0}^{+\infty} b^j.(d\omega_j - df\wedge \omega_{j+1}) \\
     & a.\sum_{j=0}^{+\infty} b^j.\omega_j  = \sum_{j=0}^{+\infty} b^j.(f.\omega_j + (j-1).\omega_{j-1}) \quad {\rm with \ the \ convention} \quad \omega_{-1} = 0 \\
     & b.\sum_{j=0}^{+\infty} b^j.\omega_j  = \sum_{j=1}^{+\infty} b^j.\omega_{j -1}
     \end{align*}
     It is easy to check that \ $D$ \ is \ $\A-$linear and that \ $D^2 = 0 $. We have a natural inclusion of complexes of left \ $\A-$modules
     $$\tilde{i} :  ({\Omega''}^{\bullet}[[b]], D^{\bullet}) \to ({\Omega'}^{\bullet}[[b]], D^{\bullet}).$$
     
     Remark that we have natural morphisms of complexes
       \begin{align*}
     & u :   (\hat{I}^{\bullet}, d^{\bullet}) \to  ({\Omega''}^{\bullet}[[b]], D^{\bullet}) \\
     & v :  (\hat{K}^{\bullet}, d^{\bullet}) \to  ({\Omega'}^{\bullet}[[b]], D^{\bullet})
    \end{align*}
    and that these morphisms are compatible with \ $i$. More precisely, this means that we have the commutative diagram of complexes
 $$   \xymatrix{ (\hat{I}^{\bullet}, d^{\bullet}) \ar[d]^i \ar[r]^u & ({\Omega''}^{\bullet}[[b]], D^{\bullet}) \ar[d]^{\tilde{i}} \\
     (\hat{K}^{\bullet}, d^{\bullet})  \ar[r]^v &  ({\Omega'}^{\bullet}[[b]], D^{\bullet}) } $$
     
     The following theorem is a variant of theorem 2.2.1. of [B. 07].

     \begin{thm}\label{(a,b)-structures}
     Let \ $X$ \ be a connected  complex manifold of dimension \ $n +1$ \ and \ $f : X \to \mathbb{C}$ \ a non constant holomorphic function with the following condition:
      $$\{ x \in X / \ df = 0 \} \subset f^{-1}(0).$$ 
      Then the morphisms of complexes \ $u$ \ and \ $v$ \ introduced above are quasi-isomorphisms. Moreover, the isomorphims that they induce on the cohomology sheaves of these complexes are compatible with the actions of \ $a$ \ and \ $b$.
  \end{thm}
  
  \smallskip
  
  This theorem builds a natural structure of left \ $\A-$modules on each of the complex \\ 
  $(\hat{K}^{\bullet}, d^{\bullet}), (\hat{I}^{\bullet}, d^{\bullet}), ([\hat{K}/\hat{I}]^{\bullet}, d^{\bullet}) $ \ and \ $f^*(\hat{\Omega}_{\mathbb{C}}^1), (\tilde{K}^{\bullet}, d^{\bullet})$ \ in the derived category of bounded complexes of sheaves of \ $\mathbb{C}-$vector spaces on \ $X$.\\
   Moreover the short exact sequences
  \begin{align*}
&  0 \to (\hat{I}^{\bullet}, d^{\bullet}) \to (\hat{K}^{\bullet}, d^{\bullet}) \to ([\hat{K}/\hat{I}]^{\bullet}, d^{\bullet}) \to 0 \tag{2}\\
&  0 \to f^*(\hat{\Omega}_{\mathbb{C}}^1) \to (\hat{K}^{\bullet}, d^{\bullet}), (\hat{I}^{\bullet}, d^{\bullet}) \to (\tilde{K}^{\bullet}, d^{\bullet}) \to 0 . \tag{3}
\end{align*}
are equivalent to short exact sequences of complexes of left \ $\A-$modules in the derived category.

  \parag{Proof} We have to prove that for any \ $p \geq 0$ \ the maps \ $\mathcal{H}^p(u)$ \ and \ $\mathcal{H}^p(v)$ \ are bijective and compatible with the actions of \ $a$ \ and \ $b$. The case of \ $\mathcal{H}^p(v)$ \ is handled (at least for \ $n \geq 2$ \ and \ $f$ \ reduced) in prop. 2.3.1. of [B.07]. To seek completness and for the convenience of the reader  we shall treat here the case of \ $\mathcal{H}^p(u)$.\\
  First we shall prove the injectivity of \ $\mathcal{H}^p(u)$. Let \ $\alpha = df\wedge \beta \in \hat{I}^p \cap Ker\, d$ \ and assume that we can find \ $U = \sum_{j=0}^{+\infty} b^j.u_j \in {\Omega''}^{p-1}[[b]]$ \ with \ $\alpha = DU$. Then we have the following relations
  $$ u_0 = df \wedge \zeta,  \quad \alpha = du_0 - df \wedge u_1 \quad {\rm and} \quad du_j = df\wedge u_{j+1} \quad \forall j \geq 1. $$
  For \ $j \geq 1$ \ we have \ $[du_j] = b[du_{j+1}]$ \ in \ $ \mathcal{H}^p(\hat{K}^{\bullet}, d^{\bullet})$;  using corollary 2.2. of [B.07] which gives the \ $b-$separation of \ $\mathcal{H}^p(\hat{K}^{\bullet}, d^{\bullet})$, this implies \ $[du_j] = 0, \forall j \geq 1$ \ in \ $ \mathcal{H}^p(\hat{K}^{\bullet}, d^{\bullet})$. For instance we can find \ $\beta_1 \in \hat{K}^{p-1}$ \ such that \ $du_1 = d\beta_1$. Now, by de Rham, we can write \ $u_1 = \beta_1 + d\xi_1$ \ for \ $p \geq 2$, where \ $\xi_1 \in \hat{\Omega}^{p-2}$. Then we conclude that
   \ $ \alpha = -df\wedge d(\xi_1 + \zeta) $ \ and \ $[\alpha] = 0$ \ in \ $\mathcal{H}^p(\hat{I}^{\bullet}, d^{\bullet})$.\\
   For \ $p = 1$ \ we have \ $u_0 = 0, \  du_1 = 0 $ \ so \ $[\alpha] = [-df\wedge d\xi_1] = 0$ \ in \ $\mathcal{H}^1(\hat{I}^{\bullet}, d^{\bullet})$.\\
   We shall show now that the image of \ $\mathcal{H}^p(u)$ \ is dense in \ $\mathcal{H}^p({\Omega''}^{\bullet}[[b]], D^{\bullet})$ \  for its  \ $b-$adic topology. Let \ $\Omega : = \sum_{j=0}^{+\infty} \ b^j.\omega_j \in {\Omega''}^{p}[[b]]$ \ such that \ $D\Omega = 0$. The following relations holds \ $ d\omega_j = df\wedge \omega_{j+1} \quad \forall j \geq 0 $ \ and \ $\omega_0 \in \hat{I}^p$. The corollary 2.2. of [B.07] again allows to find \ $\beta_j \in \hat{K}^{p-1}$ \ for any \ $j \geq 0$ \ such that \ $d\omega_j = d\beta_j$. Fix \ $N \in \mathbb{N}^*$. We have
   $$ D(\sum_{j=0}^N b^j.\omega_j) = b^N.d\omega_N = D(b^N.\beta_N) $$
   and \ $\Omega_N : = \sum_{j=0}^N b^j.\omega_j  - b^N.\beta_N $ \ is \ $D-$closed and in \ ${\Omega''}^{p}[[b]]$. As \ $\Omega - \Omega_N$ \ lies in \ $ b^N.\mathcal{H}^p({\Omega''}^{\bullet}[[b]], D^{\bullet})$, the sequence \ $(\Omega_N)_{N \geq 1}$ \ converges to \ $\Omega$ \ in \ $\mathcal{H}^p({\Omega''}^{\bullet}[[b]], D^{\bullet})$ \ for its \ $b-$adic topology. Let us show that each \ $\Omega_N$ \ is in the image of \ $\mathcal{H}^p(u)$.\\
   Write \ $\Omega_N : = \sum_{j=0}^N b^j.w_j $. The condition \ $D\Omega_N = 0$ \ implies \ $dw_N = 0$ \ and \ $dw_{N-1} = df\wedge w_N = 0$. If we write \ $w_N = dv_N$ \ we obtain \ $d(w_{N-1} + df\wedge v_N) = 0$ \ and  \ $\Omega_N - D(b^N.v_N) $ \ is of degree \ $N-1$ \ in \ $b$. For \ $N =1$ \ we are left with \ $w_0 + b.w_1 - (-df\wedge v_1 + b.dv_1) = w_0 + df\wedge v_1$ \ which is in \ $\hat{I}^p \cap Ker\, d$ \ because \ $dw_0 = df\wedge dv_1$.\\
   To conclude it is enough to know the following two facts
   \begin{enumerate}[i)]
   \item The fact  that \ $\mathcal{H}^p(\hat{I}^{\bullet}, d^{\bullet})$ \ is complete for its \ $b-$adic topology.
   \item The fact that \ $Im(\mathcal{H}^p(u)) \cap b^N.\mathcal{H}^p({\Omega''}^{\bullet}[[b]], D^{\bullet}) \subset Im(\mathcal{H}^p(u)\circ b^N) \quad \forall N \geq 1 $.
   \end{enumerate}
   Let us first conclude the proof of the surjectivity of \ $\mathcal{H}^p(u)$ \ assuming i) and ii).\\
   For any \ $[\Omega] \in \mathcal{H}^p({\Omega''}^{\bullet}[[b]], D^{\bullet})$ \ we know that there exists a sequence \ $(\alpha_N)_{N \geq 1}$ \ in \ $ \mathcal{H}^p(\hat{I}^{\bullet}, d^{\bullet})$ \ with \ $\Omega - \mathcal{H}^p(u)(\alpha_N) \in b^N.\mathcal{H}^p({\Omega''}^{\bullet}[[b]], D^{\bullet})$. Now the property ii) implies that we may choose the sequence \ $(\alpha_N)_{N \geq 1} $ \ such that \ $[\alpha_{N+1}] - [\alpha_N]$ \ lies in \ $ b^N.\mathcal{H}^p(\hat{I}^{\bullet}, d^{\bullet})$. So the property i) implies that the Cauchy sequence \ $([\alpha_N])_{N \geq 1} $ \ converges to \ $[\alpha] \in \mathcal{H}^p(\hat{I}^{\bullet}, d^{\bullet})$. Then the continuity of \ $\mathcal{H}^p(u)$ \ for the \ $b-$adic topologies coming from its \ $b-$linearity, implies \ $\mathcal{H}^p(u)([\alpha]) = [\Omega]$.\\
   The compatibility with \ $a$ \ and \ $b$ \ of the maps \ $\mathcal{H}^p(u)$ \ and \ $\mathcal{H}^p(v)$ \ is an easy exercice. 
   
   \smallskip
   
   Let us now prove properties i) and ii).\\
   The property i) is a direct consequence of the completion of \  $\mathcal{H}^p(\hat{K}^{\bullet}, d^{\bullet})$ \ for its \ $b-$adic topology given by the corollary 2.2. of [B.07] \ and the \ $b-$linear isomorphism \ $\tilde{b} $ \ between \ $\mathcal{H}^p(\hat{K}^{\bullet}, d^{\bullet})$ \ and \ $\mathcal{H}^p(\hat{I}^{\bullet}, d^{\bullet})$ \ constructed in the lemma 2.1.1. above.\\
   To prove ii) let \ $\alpha \in \hat{I}^p \cap Ker\, d $ \ and \ $N \geq 1$ \ such that
   $$ \alpha = b^N.\Omega + DU $$
   where \ $\Omega \in {\Omega''}^{p}[[b]]$ \ satisfies \ $D\Omega = 0$ \ and where \ $U \in  {\Omega''}^{p-1}[[b]]$. With obvious notations we have
   \begin{align*}
   & \alpha = du_0 -df\wedge u_1\\
   & \cdots \\
   & 0  = du_j - df\wedge u_{j+1}  \quad \forall j \in [1, N-1] \\
   & \cdots \\
   & 0 = \omega_0 + du_N - df\wedge u_{N+1} 
   \end{align*}
   which implies \ $D(u_0+ b.u_1+ \cdots + b^N.u_N) = \alpha + b^N.du_N$ \ and the fact that \ $du_N$ \ lies in \ $\hat{I}^p \cap Ker \, d$. So we conclude that \ $[\alpha] + b^N.[du_N] $ \ is in the kernel of \ $\mathcal{H}^p(u)$ \ which is \ $0$. Then \ $[\alpha] \in b^N.\mathcal{H}^p(\hat{I}^{\bullet}, d^{\bullet})$.
   $\hfill \blacksquare$
   
     \parag{Remark} The map
  $$ \beta : ({\Omega'}[[b]]^{\bullet}, D^{\bullet}) \to ({\Omega''}[[b]]^{\bullet}, D^{\bullet})$$
  defined by \ $\beta(\Omega) = b.\Omega$ \ commutes to the differentials and with the action of \ $b$. It induces the isomorphism \ $\tilde{b}$ \ of the lemma \ref{tilde b} on the cohomology sheaves. So it is a quasi-isomorphism of complexes of \ $\mathbb{C}[[b]]-$modules.\\
 To prove this fact, it is enough to verify that the diagram\\
  \xymatrix{&(\hat{K}^{\bullet}, d^{\bullet}) \ar[d]^{\tilde{b}} \ar[r]^v &  ({\Omega'}[[b]]^{\bullet}, D^{\bullet}) \ar[d]^{\beta} \\
 & (\hat{I}^{\bullet}, d^{\bullet}) \ar[r]^u &  ({\Omega''}[[b]]^{\bullet}, D^{\bullet})}\\
 \smallskip
 induces  commutative diagrams on the cohomology sheaves. \\
  But this is clear because if \ $\alpha = d\xi$ \ lies in \ $ \hat{K}^p \cap Ker \, d$ \ we have \ $D(b.\xi) = b.d\xi - df\wedge \xi $ \ so \ $\mathcal{H}^p (\beta)\circ \mathcal{H}^p (v)([\alpha]) = \mathcal{H}^p (u)\circ \mathcal{H}^p (\tilde{b})([\alpha])$ \ in \ $\mathcal{H}^p ({\Omega''}[[b]]^{\bullet}, D^{\bullet}). \hfill \blacksquare$

    \subsection{The existence theorem.} 
  
  Let us recall some basic definitions on the left modules over the algebra \ $\A$.
  
  \begin{defn}\label{(a,b)-module}
  An {\bf (a,b)-module} is a left \ $\A-$module which is free and of finite rank on the commutative sub-algebra \ $\mathbb{C}[[b]]$ \ of \ $\A$.\\
  An (a,b)-module \ $E$ \ is 
  \begin{enumerate}
  \item {\bf local} when \ $\exists N \in \mathbb{N}$ \ such that \ $a^N.E \subset b.E$;
  \item {\bf simple pole} when \ $a.E \subset b.E$;
  \item {\bf regular} when it is contained in a simple pole (a,b)-module;
  \item {\bf geometric} when it is contained in a simple pole (a,b)-module \ $E^{\sharp}$ \ such that the minimal polynomial of the action of \ $b^{-1}.a$ \ on \ $E^{\sharp}\big/ b.E^{\sharp}$ \ has its roots in \ $\mathbb{Q}^{+*}$.
  \end{enumerate}
  \end{defn}
  
  We shall give more details and examples of (a,b)-modules  in the section 3.\\
  Now let \ $E$ \ be any left \ $\A-$module, and define \ $B(E)$ \ as the \ $b-$torsion of \ $E$, that is to say
  $$ B(E) : = \{ x \in E \ / \  \exists N \quad  b^N.x = 0 \}.$$
  Define \ $A(E)$ \ as the \ $a-$torsion of \ $E$ \ and 
   $$\hat{A}(E) : = \{x \in E \ / \ \mathbb{C}[[b]].x \subset A(E) \}.$$
   Remark that \ $B(E)$ \ and \ $\hat{A}(E)$ \ are  sub-$\A-$modules of \ $E$ \ but that \ $A(E)$ \ is not stable by \ $b$. 
  
  \begin{defn}\label{petit}
  A left \ $\A-$module \ $E$ \ is {\bf small} when the following conditions hold
  \begin{enumerate}
  \item \ $E$ \ is a finite type \ $\mathbb{C}[[b]]-$module ;
  \item \ $B(E) \subset \hat{A}(E)$ ;
  \item \ $\exists N \ / \  a^N.\hat{A}(E) = 0 $ ;
  \end{enumerate}
  \end{defn}
  
  Recall that for \ $E$ \ small we have always the equality \ $ B(E) = \hat{A}(E)$ \ and that this complex vector space is finitie dimensional. The quotient \ $E/B(E)$ \ is an (a,b)-module called {\bf the associate (a,b)-module} to \ $E$.\\
  Conversely, any left \ $\A-$module \ $E$ \ such that \ $B(E)$ \ is a finite dimensional \ $\mathbb{C}-$vector space and such that \ $E/B(E)$ \ is an (a,b)-module is small.\\
  The following easy  criterium to be small will be used later :
  
  \begin{lemma}\label{crit. small}
  A left \ $\A-$module \ $E$ \ is small if and only if the following conditions hold :
  \begin{enumerate}
  \item \ $\exists N \ / \ a^N.\hat{A}(E) = 0 $ ;
  \item \ $B(E) \subset \hat{A}(E) $ ;
  \item  \ $\cap_{m\geq 0} b^m.E \subset \hat{A}(E) $ ;
  \item \ $Ker \, b$ \ and \ $Coker \, b$ \ are finite dimensional complex vector spaces.
  \end{enumerate}
  \end{lemma} 
  
  As the condition 3) in the previous lemma has been omitted in [B.07] (but this does not affect the results of this article because this lemma was used only in a case where this condition 3) was satisfied, thanks to proposition 2.2.1. of {\it loc. cit.}), we shall give the (easy) proof.
  \parag{Proof} First the conditions 1) to 4)  are obviously necessary. Conversely, assume that \ $E$ \ satisfies these four conditions. Then condition 2) implies that the action of \ $b$ \ on \ $\hat{A}(E)\big/B(E)$ \ is injective. But the condition 1) implies that \ $b^{2N} = 0$ \ on \ $\hat{A}(E) $ \ (see [B.06] ). So we conclude that \ $\hat{A}(E) = B(E) \subset Ker\, b^{2N}$ \ which is a finite dimensional complex vector space using condition 4) and an easy induction. Now \ $E/B(E)$ \ is a \ $\mathbb{C}[[b]]-$module which is separated for its \ $b-$adic topology. The finitness of \ $Coker \, b$ \ now shows that it is a free finite type \ $\mathbb{C}[[b]]-$module concluding the proof. $\hfill \blacksquare$
  
  \begin{defn}\label{geometric}
  We shall say that a left \ $\A-$module \ $E$ \ is {\bf geometric} when \ $E$ \ is small and when it associated (a,b)-module \ $E/B(E)$ \ is geometric.
  \end{defn}
  
  The main result of this section is the following theorem, which shows that the Gauss-Manin connexion of a proper holomorphic function produces geometric \ $\A-$modules associated to vanishing cycles and nearby cycles.

 \begin{thm}\label{Finitude}
 Let \ $X$ \ be a connected complex manifold of dimension \ $n + 1$ \ where \ $ n \in \mathbb{N}$, and let \ $f : X \to D$ \ be an non constant proper  holomorphic function to an open  disc \ $D$ \ in \ $\mathbb{C}$ \ with center \ $0$. Let us assume that \ $df$ \ is nowhere vanishing outside of \ $X_0 : = f^{-1}(0)$.\\
 Then the \ $\A-$modules 
 $$ \mathbb{H}^j(X, (\hat{K}^{\bullet}, d^{\bullet})) \quad {\rm and} \quad \mathbb{H}^j(X, (\hat{I}^{\bullet}, d^{\bullet})) $$
 are geometric for any \ $j \geq 0 $.
 \end{thm}
 
In the proof we shall use the \ $\mathcal{C}^{\infty}$ \ version of the complex \ $(\hat{K}^{\bullet}, d^{\bullet})$. We define \ $K_{\infty}^p$ \ as the kernel of \ $\wedge df : \mathcal{C}^{\infty,p} \to \mathcal{C}^{\infty,p+1}$ \ where \ $\mathcal{C}^{\infty,j}$ \ denote the sheaf of \ $\mathcal{C}^{\infty}-$ \ forms on \ $X$ \ of degree j, let \ $\hat{K}_{\infty}^p$ \ be its formal \ $f-$completion and \ $(\hat{K}_{\infty}^{\bullet}, d^{\bullet})$ \ the corresponding de Rham complex.

  The next lemma is proved in [B.07] (lemma  6.1.1.)
 
 \begin{lemma}\label{diff}
 The natural inclusion
 $$ (\hat{K}^{\bullet}, d^{\bullet}) \hookrightarrow (\hat{K}_{\infty}^{\bullet}, d^{\bullet}) $$
is a quasi-isomorphism.
 \end{lemma}
 
 \parag{Remark} As the sheaves \ $\hat{K}_{\infty}^{\bullet}$ \ are fine, so we have a natural isomorphism
 $$ \mathbb{H}^p(X, (\hat{K}^{\bullet}, d^{\bullet})) \simeq H^p\big(\Gamma(X, \hat{K}_{\infty}^{\bullet}), d^{\bullet}\big).$$
 
 Let us denote by \ $X_1$ \ the generic fiber of \ $f$. Then \ $X_1$ \ is a smooth compact complex manifold of dimension \ $n$ \ and the restriction of \ $f$ \ to \ $f^{-1}(D^*)$ \ is a locally trivial \ $\mathcal{C}^{\infty}$ \ bundle with typical fiber \ $X_1$ \ on \ $D^* = D \setminus \{0\}$, if the disc \ $D$ \ is small enough around \ $0$. Fix now \ $\gamma \in H_p(X_1, \mathbb{C})$ \ and let \ $(\gamma_s)_{s \in D^*}$ \ the corresponding multivalued horizontal family of \ $p-$cycles \ $\gamma_s \in H_p(X_s, \mathbb{C})$. Then for \ $\omega \in \Gamma(X, \hat{K}_{\infty}^p \cap Ker\, d) $ \ define the multivalued holomorphic function
 $$ F_{\omega}(s) : = \int_{\gamma_s} \frac{\omega}{df} .$$
 Let now 
  $$\Xi : = \underset{\alpha \in \mathbb{Q}\, \cap ]-1,0], \, j \in [0,n]}{\sum} \quad  \mathbb{C}[[s]].s^{\alpha}.\frac{(Log s)^j}{j!} .$$
  This is an \ $\A-$modules with \ $a$ \ acting as multiplication by \ $s$ \ and \ $b$ \ as the primitive in \ $s$ \ without constant. Now if \ $\hat{F}_{\omega}$ \ is the asymptotic expansion at \ $0$ \ of \ $F_{\omega}$, it is an element  in \ $\Xi$, and we obtain in this way an \ $\A-$linear map
 $$ Int :   \mathbb{H}^p(X, (\hat{K}^{\bullet}, d^{\bullet})) \to H^p(X_1, \mathbb{C}) \otimes_{\mathbb{C}} \Xi .$$
 To simplify notations, let  \ $E : =  \mathbb{H}^p(X, (\hat{K}^{\bullet}, d^{\bullet}))$. Now using Grothendieck theorem [G.66], there exists \ $N \in \mathbb{N}$ \ such that \ $Int(\omega) \equiv 0 $, implies \ $a^N.[\omega] = 0$ \ in \ $E$.  As the converse is clear we conclude that \ $\hat{A}(E) =  Ker(Int)$. It is also clear that \ $B(E) \subset Ker(Int)$ \ because \ $\Xi$ \ has no \ $b-$torsion. So we conclude that \ $E$ \ satisfies properties 1) and 2) of the lemma \ref{crit. small}.\\
 The property 3) is also true because of the regularity of the Gauss-Manin connexion of \ $f$.
 
 \parag{End of the proof of theorem \ref{Finitude}} To show that \ $E : =  \mathbb{H}^p(X, (\hat{K}^{\bullet}, d^{\bullet}))$ \ is small, it is enough to prove that \ $E$ \ satisfies the condition 4) of the lemma \ref{crit. small}. Consider now the long exact sequence of hypercohomology of the exact sequence of complexes
 $$   0 \to (\hat{I}^{\bullet}, d^{\bullet}) \to (\hat{K}^{\bullet}, d^{\bullet}) \to ([\hat{K}/\hat{I}]^{\bullet}, d^{\bullet}) \to 0 .$$
 It contains the exact sequence
 $$  \mathbb{H}^{p-1}(X, ([\hat{K}\big/\hat{I}]^{\bullet}, d^{\bullet})) \to \mathbb{H}^p(X, (\hat{I}^{\bullet}, d^{\bullet})) \overset{\mathbb{H}^p(i)}{\to} \mathbb{H}^p(X, (\hat{K}^{\bullet}, d^{\bullet})) \to \mathbb{H}^{p}(X, ([\hat{K}\big/\hat{I}]^{\bullet}, d^{\bullet})) $$
 and we know that \ $b$ \ is induced on the complex of \ $\A-$modules  quasi-isomorphic to \ $(\hat{K}^{\bullet}, d^{\bullet})$ \ by the composition \ $i\circ \tilde{b}$ \ where \ $\tilde{b}$ \ is a quasi-isomorphism of complexes of \ $\mathbb{C}[[b]]-$modules. This implies that the kernel  and the cokernel of \ $\mathbb{H}^p(i)$ \ are isomorphic (as \ $\mathbb{C}-$vector spaces) to \ $Ker\, b$ \ and \ $Coker \, b$ \ respectively. Now to prove that \ $E$ \ satisfies condition 4) of the lemma \ref{crit. small} it is enough to prove finite dimensionality for the vector spaces \ $  \mathbb{H}^{j}(X, ([\hat{K}\big/\hat{I}]^{\bullet}, d^{\bullet})) $ \ for all \ $j \geq 0 $.\\
 But the sheaves \ $[\hat{K}\big/\hat{I}]^j \simeq [Ker\,df\big/Im\, df]^j$ \ are coherent on \ $X$ \ and supported in \ $X_0$. The spectral sequence
 $$ E_2^{p,q} : = H^q\big( H^p(X, [\hat{K}\big/\hat{I}]^{\bullet}), d^{\bullet}\big) $$ 
 which converges to \ $ \mathbb{H}^{j}(X, ([\hat{K}\big/\hat{I}]^{\bullet}, d^{\bullet})) $, is a bounded complex of finite dimensional vector spaces by Cartan-Serre. This gives the desired finite dimensionality.\\
 To conclude the proof, we want to show that \ $E/B(E)$ \ is geometric. But this is an easy consequence of the regularity of the Gauss-Manin connexion of \ $f$ \ and of the Monodromy theorem, which are already incoded in the definition of \ $\Xi$ : the injectivity on \ $E/B(E)$ \ of the \ $\A-linear$ \ map \ $Int$ \ implies that \ $E/B(E)$ \ is geometric. \\
 Remark now that the piece of  exact sequence above gives also the fact that \ $\mathbb{H}^p(X, (\hat{I}^{\bullet}, d^{\bullet}))$ \ is geometric, because it is an exact sequence of \ $\A-$modules. $\hfill \blacksquare$
 
 \section{Basic properties.}

\subsection{Definition and examples.}

First recall in a more na{\"i}ve way  the definition of an (a,b)-module.

\begin{defn}\label{(a,b)_module}
An {\bf (a,b)-module \ $E$} \ is a free finite type \ $\mathbb{C}[[b]]-$module with a \ $\mathbb{C}-$linear endomorphism \ $a : E \to E$ \ which is continuous for the \ $b-$adic topology of \ $E$ \ and satisfies \ $a.b - b.a = b^2$.\\
The {\bf rank of \ $E$}, denote by \ $rank(E)$,  will be the rank of \ $E$ \ as a \ $\mathbb{C}[[b]]-$module.
\end{defn}

\parag{Remarks}
\begin{enumerate}
\item Let \ $(e_1, \cdots, e_k)$ \ a \ $\mathbb{C}[[b]]-$basis of a free finite type  $\mathbb{C}[[b]]-$module. Then choosing arbitrarily elements \ $(\varepsilon_1, \cdots, \varepsilon_k)$ \ and defining \ $a.e_j = \varepsilon_j \quad \forall j \in [1,k]$ \ gives an (a,b)-module: the commutation relation implies that \ $\forall n \in \mathbb{N}$ \ we have \ $ a.b^n = b^n.a + n.b^{n+1} $ \ so \ $a$ \ is defined on \ $\sum_{j =1}^k \quad \mathbb{C}[b].e_j$. The continuity assumption gives its (unique) extension.
\item There is a natural (a,b)-module associated to every algebraic linear differential system (see [B.95] p.42)
$$ Q(z).\frac{dF}{dz} = M(z).F(z), \quad Q \in \mathbb{C}[z], \quad M \in End(\mathbb{C}^n)\otimes_{\mathbb{C}} \mathbb{C}[z] .$$
\end{enumerate}

In the sequel of this article we shall mainly consider regular (a,b)-modules (see definition recalled below).  To try to convince the reader that the ''general'' (a,b)-module structure is interesting, let me quote the following result, which is quite elementary in the regular case, but which is not so easy in general.

\begin{thm}([B.95] th.1bis p.31)
Let \ $E$ \ be an (a,b)-module. Then the kernel and cokernel of ''a''  are finite dimensional.
\end{thm}

This result implies a general finiteness theorem for extensions of (a,b)-modules (see [B.95] and also section 1.3).

\bigskip

\begin{defn}\label{simple pole}
We shall say that an (a,b)-module \ $E$ \ has a {\bf simple pole} when the inclusion \ $a.E \subset b.E$ \ is satisfied.
\end{defn}

This terminology comes from the terminology of meromorphic connexions (see for instance [D.70]). 

\parag{Example}
 For any \ $\lambda \in \mathbb{C}$ \ define the simple pole  rank 1 \ (a,b)-module \ $E_{\lambda}$ \ as \ $E : = \mathbb{C}[[b]].e_{\lambda}$ \ where ''$a$'' is defined by the relation \ $a.e_{\lambda} = \lambda.b.e_{\lambda}$. $\hfill \square$
 
 \bigskip
 
 As an introduction to our second  theorem, the reader may solve the following exercice by direct computation.
 \parag{Exercice}
For any \ $S \in \mathbb{C}[[b]]$ \ show that the simple pole (a,b)-module defined by  \ $E : = \mathbb{C}[[b]].e_S $ \ and \ $a.e_S = b.S(b).e_S$ \ is isomorphic to \ $E_{\lambda}$ \ with \ $\lambda = S(0)$ \\
(hint: begin by looking for \ $ \alpha_1 \in \mathbb{C}$ \ such that \ $(a - S(0).b)(e + \alpha_1.b.e) \in b^3.E$). $\hfill \square$

\bigskip

For a simple pole (a,b)-module, the linear map \ $b^{-1}.a : E \to E $ \ is well defined and induces an endomorphism \ $ f : = b^{-1}.a : E/b.E \to E/b.E$. For any \ $\lambda \in \mathbb{C}$ \ we shall denote by \ $\lambda_{min}$ \ the smallest eigenvalue of \ $f$ \ which is in \ $\lambda + \mathbb{Z}$. Then for \ $\lambda = \lambda_{min}- k$ \ with \ $k \in \mathbb{N}^*$ \ the bijectivity of the map \ $f - \lambda$ \ on \ $E/b.E$ \ implies easily its bijectivity on \ $E$ \ (see the exercice above). It gives then the equality 
$$ (a - \lambda.b).E = b.E.$$
Using this remark, it is not difficult to prove the following result from [B.93] (prop.1.3. p.11) that we shall use  later on.

\begin{prop}\label{sub min}
Let \ $E$ \ be a simple pole (a,b)-module, and let \ $\lambda \in \mathbb{C}$ \ and \ $\kappa \in \mathbb{N}$ \ such that \ $\lambda - \kappa \leq \lambda_{min} $. If \ $y \in E$ \ satisfies \ $(a - \lambda.b).y \in b^{\kappa+2}.E $ \ then there exists an unique \ $\tilde{y} \in E$ \ such that \ $(a - \lambda.b).\tilde{y} = 0 $ \ and \ $\tilde{y} - y \in b^{\kappa+1}.E$.
\end{prop}

An easy consequence of this proposition is that for an eigenvalue \ $\lambda$ \ of \ $f$ \ such that \ $\lambda = \lambda_{min}$ \ there always exists a non zero \ $x \in E$ \ such that \ $(a- \lambda.b).x = 0$. This gives an embedding of \ $E_{\lambda}$ \ in \ $E$. Remark also that if \ $E$ \ is a non zero simple pole (a,b)-module, such a \ $\lambda$ \ always exists. This leads to a rather precise description a of ''general'' simple pole (a,b)-module (see [B.93] th. 1.1 p.15). 

\bigskip

\begin{defn}\label{regular}
An (a.b)-module \ $E$ \ is  {\bf regular} \ when its saturation by \ $b^{-1}.a$ \ in \ $E[b^{-1}]$ \ is finitely generated on \ $\mathbb{C}[[b]]$.
\end{defn}

\bigskip

We shall denote \ $E^{\sharp}$ \ this saturation. It is a simple pole (a,b)-module and it is the smallest simple pole (a,b)-module containing \ $E$ \ in the sense that for any  (a,b)-linear morphism  \ $j : E \to F$ \ where \ $F$ \ is a simple pole (a,b)-module, there exists a unique (a,b)-linear extension \ $j^{\sharp} : E^{\sharp} \to F$ \ of \ $j$. \\

\smallskip

It is easy  to show that a regular (a,b)-module of rank 1  is isomorphic to some \ $E_{\lambda}$ \ for some \ $\lambda \in \mathbb{C}$. The classification of rank 2 regular (a,b)-module is not so obvious. We recall it here for a later use

\begin{prop}(see [B.93] prop.2.4 p. 34)\label{class}
The list of rank 2 regular (a,b)-modules is, up to isomorphism, the following :
\begin{enumerate}
\item \ $E_{\lambda} \oplus E_{\mu} $ \ for \ $(\lambda,\mu) \in \mathbb{C}^2/\frak{S}_2 $.
\item For any \ $\lambda \in \mathbb{C}$ \ and any \ $n \in \mathbb{N}$ \ let \ $E_{\lambda}(n)$ \ be the simple pole (a,b)-module with basis \ $(x,y)$ \ such that
$$  a.x = (\lambda + n).b.x + b^{n+1}.y \quad {\rm and} \quad  a.y = \lambda.b.y .$$
\item For any \ $(\lambda,\mu) \in \mathbb{C}^2/\frak{S}_2 $ \ let \ $E_{\lambda,\mu}$ \ the rank 2 regular (a,b)-module with basis \ $(y,t)$ \ such that
$$ a.y = \mu.b.y \quad {\rm and } \quad a.t = y + (\lambda-1).b.t .$$
\item For any \ $\lambda \in \mathbb{C}$, any \ $n \in \mathbb{N}^*$ \ and any \ $\alpha \in \mathbb{C}^*$ \ let \ $E_{\lambda, \lambda-n}(\alpha)$ \ be the rank 2 regular (a,b)-module with basis \ $(y,t)$ \ such that
$$ a.y = (\lambda-n).b.y \quad {\rm and} \quad a.t = y + (\lambda-1)b.t + \alpha.b^n.y $$
\end{enumerate}
\end{prop}

Note that the first two cases are simple pole (a,b)-modules.\\
 The saturation by \ $b^{-1}.a$ \  in case 3.  is generated by \ $b^{-1}.y$ \ and \ $t$ \ as a \ $\mathbb{C}[[b]]-$module. It is isomorphic to \ $E_{\lambda-1} \oplus E_{\mu-1}$ \ for \ $\lambda \not= \mu$ \ and to \ $E_{\lambda-1}(0)$ \ for \ $\lambda = \mu$.\\
 The saturation  by \ $b^{-1}.a$ \  in case 4.  is generated by \ $b^{-1}.y$ \ and \ $t$ \ as a \ $\mathbb{C}[[b]]-$module. It is isomorphic to \ $E_{\lambda-n-1}(n)$ \ for any non zero value of \ $\alpha$.
 
 \bigskip

To conclude this first section, let me recall also the theorem of existence of Jordan-H{\"o}lder sequences for regular (a,b)-module, which will be usefull in the induction in the proof of  our  result .

\begin{thm}(see [B.93] th. 2.1 p.30)\label{J-H}
For any regular  rank k (a,b)-module \ $E$ \ there exists a sequence of sub-(a,b)-modules
$$ 0 = E^0 \subset E^1 \subset \cdots \subset E^{k-1} \subset E^k = E $$
such that for any \ $j \in [1,k]$ \ the quotient \ $E^j/E^{j-1}$ \ is isomorphic to \ $E_{\lambda_j}$. Moreover we may choose for \ $E^1$ \ any normal\footnote{normal means \ $E^1 \cap b.E = b.E^1$, so that \ $E/E^1$ \ is again free on \ $\mathbb{C}[[b]]$.}  rank 1 sub-(a,b)-module of \ $E$.\\
The number \ $\alpha(E) : = \sum_{j=1}^k \ \lambda_j $ \ is independant of the choice of the Jordan-H{\"o}lder sequence. It is given by the following formula
$$ \alpha(E) = trace\big(b^{-1}.a : E^{\sharp}/b.E^{\sharp} \to E^{\sharp}/b.E^{\sharp} \big) + \dim_{\mathbb{C}}(E^{\sharp}/E) .$$
\end{thm}

\subsection{The regularity order.}

\begin{defn}\label{order of regularity}
Let \ $E$ \ be a regular (a,b)-module. We define the {\bf regularity order of \ $E$} as the smallest integer \ $k \in \mathbb{N}$ \ such that the inclusion
\begin{equation*}
 a^{k+1}.E \subset \sum_{j=0}^k \quad a^j.b^{k-j+1}.E \tag{reg.}
 \end{equation*}
is valid. We shall note this integer \ $or(E)$.\\
We define also {\bf the index \ $\delta(E)$ \ of \ $E$} \ as the smallest integer \ $m \in \mathbb{N}$ \ such that \ $E^{\sharp} \subset b^{-m}.E $.
\end{defn}

\parag{Remarks}
\begin{enumerate}[i)]
\item The (a,b)-module \ $E$ \ has a simple pole if an only iff \ $or(E) = 0 $.
\item The inclusion \ (reg.) \ implies that \ $ (b^{-1}.a)^{k+1}.E \subset \Phi_k(E) : = \sum_{j=0}^k \quad (b^{-1}.a)^j.E$ \ and this implies that \ $\Phi_k(E)$ \ is stable by \ $b^{-1}.a$. So \ $\Phi_k(E)$ \ is a simple pole (a,b)-module contained in \ $b^{-k}.E \subset E[b^{-1}]$. This implies clearly the regularity of \ $E$. \\
For \ $k = or(E)$ \ we have \ $E^{\sharp} = \Phi_k(E) \subset b^{-k}.E$. So we have \ $\delta(E) \leq or(E)$.
\item As the quotient \ $b^{-k}.E/E$ \ is a finite dimensional \ $\mathbb{C}-$vector space, the quotient \ $E^{\sharp}/E$ \ is always a finite dimensional \ $\mathbb{C}-$vector space. $\hfill \square$
\end{enumerate}

The remark iii) shows that for a regular (a,b)-module \ $E$ \  there always exists a simple pole sub-(a,b)-module of \ $E$ \ which is  a finite codimensional vector space in \ $E$. This comes from the fact that for \ $k =  \delta(E)$ \ we have \ $ b^k.E^{\sharp} \subset E $ \ and that  \ $b^k.E^{\sharp}$ \ has a simple pole.

\parag{Example} The inequality \ $\delta(E) \leq or(E)$ \ may be strict for \ $or(E) \geq 2$. For instance the (a,b)-module of rank 3  with \ $\mathbb{C}[[b]]-$basis \ $e_1, e_2, e_3$ \ with \\
 $a.e_1 = e_2, \quad a.e_2 = b.e_3, \quad a.e_3 = 0 $ \ has index 1  and regularity order 2 : an easy computation gives that  a \ $\mathbb{C}[[b]]-$basis for \ $E^{\sharp} $ \ is given by \ $e_1, b^{-1}.e_2, b^{-1}.e_3$, and that a  \ $\mathbb{C}[[b]]-$basis for \ $E + b^{-1}.a.E$ \ is given by \ $e_1, b^{-1}.e_2, e_3$. $\hfill \square$

\begin{defn}\label{$E^b$}
Let \ $E$ \ be a regular (a,b)-module. The {\bf biggest simple pole sub-(a,b)-module of \ $E$} \ exists and has finite \ $\mathbb{C}-$codimension in \ $E$. We shall note it \ $E^b$.
\end{defn}

In general, for \ $k = \delta(E)$ \  the inclusion \ $b^k.E^{\sharp} \subset E^b $ \ is strict. For instance this is the case for \ $E_{\lambda,\mu} \oplus E_{\nu}$.

\bigskip

\begin{lemma}\label{inclusion} 
Let \ $E$ \ be a regular (a,b)-module. The smallest integer \ $m$ \ such we have \ $b^m.E \subset E^b$ \ is equal to \ $\delta(E)$.
\end{lemma}

\parag{Proof} Let \ $k : = \delta(E)$. Then \ $b^k.E^{\sharp}$ \ is a simple pole sub-(a,b)-module of \ $E$. So we have \ $b^k.E \subset b^k.E^{\sharp} \subset E^b$. Conversely, an inclusion \ $b^m.E \subset E^b $ \ gives \ $E \subset b^{-m}.E^b $. As \ $b^{-m}.E^b$ \ has a simple pole this implies \ $E^{\sharp} \subset b^{-m}.E^b \subset b^{-m}.E $. So \ $\delta(E) \leq m$.  $\hfill \blacksquare$

\parag{Examples} In the case 3 of the proposition \ref{class} \ $E^b$ \ is generated as a \ $\mathbb{C}[[b]]-$module by \ $y$ \ and \ $b.t$, so \ $E^b = b.E^{\sharp}$.\\
In case 4 we have also \ $E^b = b.E^{\sharp}$.

\bigskip

\begin{lemma}\label{ordre reg.}
Let \ $E$ \ be a regular (a,b)-module. For any exact sequence of (a,b)-modules
\begin{equation*}
 0 \to E' \to E \overset{\pi}{ \to} E''  \to 0 \tag{*}
\end{equation*}
we have  \ $or(E'') \leq or(E) \leq rank(E') +  or(E'') $.\\
As a consequence, the order of regularity of \ $E$ \ is at most \ $rank(E) - 1$ \ for any regular non zero (a,b)-module. 
\end{lemma}

\parag{Proof}  The inequality \ $or(E'') \leq or(E)$ \ is trivial because an inequality
$$  a^{k+1}.E   \subset \sum_{j=0}^{k} \ a^j.b^{k-j+1}.E $$
implies the same for \ $E''$ \ and, by definition, the best such integer \ $k$ \ is the order of regularity.\\
The crucial case is when \ $E'$ \ is of rank 1 . So we may assume that \ $E' \simeq E_{\lambda}$ \ for some \ $\lambda \in \mathbb{C}$ \ (see \ref{J-H} or  [B.93] prop.2.2 p.23). Let \ $k = or(E'')$. Then the inclusion 
\begin{equation*}
 a^{k+1}.E'' \subset \sum_{j=0}^{k} \ a^j.b^{k-j+1}.E'' \tag{1}
 \end{equation*}
 implies that 
 \begin{equation*}
 a^{k+1}.E \subset \sum_{j=0}^{k} \ a^j.b^{k-j+1}.E + b^l.E_{\lambda}  \tag{2}
 \end{equation*}
  for some \ $l \in \mathbb{N}$. In fact we can take for \ $l$ \  the smallest integer such that the generator \ $e_{\lambda}$ \ of \ $E_{\lambda}$ \ (defined up to \ $\mathbb{C}^*$ \ by the relation \ $a.e_{\lambda} = \lambda.b.e_{\lambda}$) satisfies \ $ b^l.e_{\lambda} \in \Psi_k  = \sum_{j=0}^{k} \ a^j.b^{k-j+1}.E $. \\ 
  Remark that this integer \ $l \geq 0$ \ is well defined because \ $b^{k+1}.e_{\lambda} \in \Psi_k$. Moreover, as \ $\Psi_k$ \ is a \ $\mathbb{C}[[b]]-$submodule of \ $E$, \ $b^l.e_{\lambda} \in \Psi_k$  \ implies \ $b^l.E_{\lambda} \subset \Psi_k$.\\
 Now, thanks to \ $(2)$ \  we have 
 \begin{equation*}
 a^{k+2}.E \subset \sum_{j=0}^k \ a^{j+1}.b^{k+1-j}.E \quad + a.b^l.E_{\lambda} \tag{3}
 \end{equation*}
 which gives
 \begin{equation*}
 a^{k+2}.E \subset \sum_{j=0}^{k+1} \ a^j.b^{k-j+2}.E  \tag{4}
 \end{equation*}
 because \ $a.b^l.E_{\lambda} = b.b^l.E_{\lambda} \subset b.\Psi_k $.\\
 This proves that \ $or(E) $ \ is at most \ $k + 1 = or(E'') + rank(E') $.\\
 Assume now that our inequality is proved for \ $E'$ \ of rank \ $p - 1$ \ and consider an exact sequence \ $(^*)$ \ with  \ $rank(E')$ \ equal \ $p \geq 2$. Let \ $E_{\lambda} \subset E'$ \ be a normal rank 1  sub-(a,b)-module of \ $E'$ \ (see \ref{J-H} or  [B.93]  prop.2.2 p.23 for a proof of the existence of such sub-(a,b)-module) and consider the exact sequence of (a,b)-modules (using the fact that \ $E_{\lambda}$ \ is also normal in \ $E$; see lemma 2.5 of [B.93])
 \begin{equation*}
 0 \to E'/E_{\lambda} \to E/E_{\lambda} \to E'' \to 0 
 \end{equation*}
 Using the induction hypothesis and the rank 1 case we get
 $$ or(E) \leq   or(E/E_{\lambda}) + 1 \leq p-1 + or(E'') + 1 = p + or(E'') .$$
 Now using an easy induction  (or a Jordan-H{\"o}lder sequence for \ $E$)  we obtain \ $or(E) \leq rank(E) - 1$ \ for any regular \ $E$. $\hfill \blacksquare$
 
 \parag{Remark} In the situation of the previous lemma we have \ $\delta(E') \leq \delta(E)$. This is a consequence of the obvious inclusion \ $ (E')^{\sharp} \subset E'[b^{-1}] \cap E^{\sharp}$ : assume that \ $x \in E'[b^{-1}] \cap E^{\sharp}$ ; then, for \ $k : = \delta(E)$ \ we have  \ $ b^k.x\in E'[b^{-1}] \cap E$ \ so that \ $b^{N+k}.x \in E'$ \ for \ $N$ \ large enough. As \ $E/E'$ \ has no $b-$torsion, we conclude that \ $b^k.x \in E'$. So our initial inclusion implies \ $\delta(E') \leq k$. $\hfill \square$
 
 \bigskip
 
 \subsection{Duality.}
 
 In this section we consider more carefully  the associative and unitary \ $\mathbb{C}-$algebra
 $$ \tilde{\mathcal{A}} : = \big\{ \sum_0^{\infty} \  P_n(a).b^n \quad {\rm with} \quad P_n \in \mathbb{C}[z] \big\} $$
with the commutation relation \ $a.b - b.a = b^2 $, and such that the left and right multiplications by \ $a$ \ are continuous for the \ $b-$adic topology\footnote{remark that for each \ $k \in \mathbb{N}$ \ $b^k.\A = \A.b^k$.} of \ $\A$.

\parag{The right structure as a commuting left-structure on \ $\A$} \quad \\
There exits an unique \ $\mathbb{C}-$linear (bijective) map \ $\theta : \A \to \A$ \ with the following properties
\begin{enumerate}[i)]
\item \ $\theta(1) = 1, \quad \ \theta(a) = a, \quad \theta(b) = - b $;
\item \ $\theta(x.y) = \theta(y).\theta(x) \quad \forall x,y \in \A $.
\item \ $\theta$ \ is continuous for the \ $b-$adic topology of \ $\A$
\end{enumerate}
The uniqueness is an easy consequence of iii) and  the fact that the conditions i) and ii) implies \ $\theta(b^p.a^q) = (-1)^p.a^q.b^p \quad \forall p,q \in \mathbb{N}$. Existence is then clear from the explicit formula deduced from this remark.\\
We  define a new structure of left \ $\A-$module on \ $\A$, {\bf called the \ $\theta-$structure} and denote by \ $x_*\square$,  by the formula
$$ x_*y = y.\theta(x) .$$
It is easy to see that this new left-structure on \ $\A$ \ commutes with the ordinary one and that with this \ $\theta-$structure \ $\A$ \ is still free of rank one as a left \ $\A-$module.

\begin{defn}\label{Hom}
Let \ $E$ \ be a (left) \ $\A-$module. On the \ $\mathbb{C}-$vector space \ $Hom_{\A}(E,\A)$ \ we define a  left  \ $\A-$module structure using the \ $\theta-$structure on \ $\A$. Explicitely this means that for \ $\varphi \in Hom_{\A}(E,\A)$ \ and \ $x \in \A$ \ we let
$$\forall e \in E \quad  (x.\varphi)(e) : = x_*\varphi(e) = \varphi(e).\theta(x).$$
We obtain in this way a left  \ $\A-$module that we shall still denote \ $Hom_{\A}(E,\A)$.
\end{defn}

It is clear that \ $E \to Hom_{\A}(E,\A)$ \ is a contravariant functor which is left exact in the category of left \ $\A-$modules. As every finite type left\ $\A-$module has a resolution of length \ $\leq 2$ \ by free finite type modules ( see [B.95] cor.2 p.29), we shall denote by \ $Ext^i_{\A}(E, \A), i \in [0,2]$ \ the right derived functors of this functor. They are finite type left \ $\A-$modules when \ $E$ \ is finitely generated because \ $\A$ \ is left noetherian (see [B.95] prop.2 p.26).

\bigskip

Any (a,b)-module is a left \ $\A-$module. They are characterized by the existence of special simple resolutions.

 \begin{lemma}\label{Resol.}
 Let \ $M$ \ be a \  $(p,p)$ \ matrix with entries in the ring \ $\mathbb{C}[[b]] $. Then the left \ $\tilde{\mathcal{A}}-$linear map \ $ Id_p.a - M : \tilde{\mathcal{A}}^p \to \tilde{\mathcal{A}}^p$ \ given by \\ $$^tX : = (x_1, \cdots,x_p) \  \to  \  ^tX.(Id_p.a - M)$$
 is injective. Its cokernel is the (a,b)-module \ $E$ \ given as follows : \\ 
 $E$ \ has a \ $\mathbb{C}[[b]]$ \ base \ $e : = (e_1, \cdots, e_p)$ \ and \ $a$ \ is defined by the two conditions
 \begin{enumerate}
 \item \ $a.e : = M(b).e $ ;
 \item the left action of \ $a$ \ is continuous for the \ $b-$adic topology of \ $E$.
 \end{enumerate}
 Any (a,b)-module is obtained in this way and so, as a \ $\A-$left-module, has a resolution of the form
 \begin{equation*}
  0 \to \A^p \  \overset{^t\square.(Id_p.a - M)}{\longrightarrow} \ \A^p \to E \to 0 . \tag{@}
  \end{equation*}
  \end{lemma}
  
  \parag{Proof} First remark that for \ $x \in \A$ \ the condition \ $x.a \in b.\A $ \ implies \ $x \in b.\A $. Now let us prove, by induction on \ $n \geq 1 $, that, for any \ $(p,p)$ \ matrix \ $M$ \ with entries in \ $\mathbb{C}[[b]]$ \ the condition \ $ ^tX.(Id_p.a - M) = 0 $ \ implies \ $^tX \in b^n.\A^p$.\\
   For \ $n =1$ \ this comes from the previous remark. Let assume that the assertion is proved for \ $n \geq 1$ \ and consider an \ $X \in \A^p$ \ such that \ $ ^tX.(Id_p.a - M) = 0 $. Using the induction hypothesis we can find \ $Y \in \A^p$ \ such that \ $X = b^n.Y$. Now we obtain, using \ $a.b^n = b^n.a + n.b^{n+1}$ \  and the fact that \ $\A$ \ has no zero divisor, the relation
   $$ ^tY(Id_p.a -(M + n.Id_p.b)) = 0 $$
   and using again our initial remark we conclude that \ $Y \in b.\A^p $ \ so \ $X \in b^{n+1}.\A^p$.\\
   So such an \ $X$ \ is in \ $\cap_{n\geq 1} \ b^n.\A^p = (0)$. \\
   The other assertions of the lemma are obvious. $\hfill \blacksquare$
   
   \bigskip
   
   We recall now a construction given in [B.95] which allows to compute more easily the vector spaces \ $Ext^i_{\A}(E,F)$ \ when \ $E,F$ \ are (a,b)-modules 
   
   \begin{defn}\label{Hom(a,b)}
   Let \ $E,F$ \ two (a,b)-modules. Then the \ $\mathbb{C}[[b]]-$module \ $Hom_b(E,F)$ \ is again a free and finitely generated \ $\mathbb{C}[[b]]-$module. Define on it an (a,b)-module structure in the following way.
   \begin{enumerate}
   \item First change the sign of the action of \ $b$. So \ $S(b) \in \mathbb{C}[[b]]$ \ will act as \ $\check{S}(b) = S(-b)$.
   \item Define \ $a$ \ using the linear map
   $ \Lambda : Hom_b(E,F) \to Hom_b(E,F)$ \ given by \ $\Lambda(\varphi)(e) = \varphi(a.e) - a.\varphi(e)$.
   \end{enumerate}
   We shall denote \ $Hom_{a,b}(E,F)$ \ the corresponding (a,b)-module.
   \end{defn}
   The verification that \ $\Lambda(\varphi)$ \ is \ $\mathbb{C}[[b]]-$linear and that \ $\Lambda.\check{b} - \check{b}.\Lambda = \check{b}^2$ \ are easy (and may be found in [B.95] p.31).
   
  \parag{Remark} In {\it loc. cit.} we defined the (a,b)-module structure on \ $Hom_{a,b}(E,F)$ \ with opposite  signs for \ $a$ \ and \ $b$. The present convention is better because it fits with the usual definition of the formal adjoint of a differential operator : \ $z^* = z$ \ and \ $(\partial/\partial z)^* = -\partial/\partial z$. $\hfill \square$
  
  \bigskip

   The following lemma is also proved in {\it loc.cit.}
   
   \begin{lemma}\label{hom}
    Let \ $E,F$ \ two (a,b)-modules. Then there is a functorial isomorphism of \ $\mathbb{C}-$vector spaces 
     $$ H^i\Big( Hom_{a,b}(E,F) \overset{a}{\to}Hom_{a,b}(E,F) \Big) \to Ext^i_{\A}(E, F)\quad \forall i \geq 0 .$$
     Here the map \ $a$ \ of the complex \ $Hom_{a,b}(E,F) \overset{a}{\to}Hom_{a,b}(E,F) )$ \ is equal to the \ $\Lambda$ \ defined above which is, by definition,  the operator \ $''a''$ \ of the (a,b)-module \ $Hom_{a,b}(E,F)$.
     \end{lemma}
     
     Now the following corollary of the lemma  \ref{Resol.} gives that the two natural ways of defining the dual of an (a,b)-module give the same answer.  
   
     \begin{cor}\label{Dual Reg.}
   Let \ $E$ \ an (a,b)-module. There is a functorial isomorphism of (a,b)-modules between the following two (a,b)-modules constructed as follows : 
   \begin{enumerate}
   \item \ $Ext^1_{\A}(E,\A)$ \ with the \ $\A-$structure defined by the \ $\theta-$structure of \ $\A$.
   \item \ ${\it Hom}_{a,b}(E, E_0) $ \ where \ $E_0 : = \A \big/ \A.a$.
   \end{enumerate}
   \end{cor}
   
   \parag{Proof} Using a free resolution \ $(@)$ \ of \ $E$ \ deduced from a \ $\mathbb{C}[[b]]-$basis \\
    $e : = (e_1, \cdots, e_p)$ \ we obtain, by the previous lemma, an exact sequence
    \begin{equation*}
  0 \to \A^p \  \overset{(Id_p.a - ^tM).\square}{\longrightarrow} \ \A^p \to Ext^1_{\A}(E, \A) \to 0 . \tag{@@}
  \end{equation*}
  of left \ $\A-$modules where \ $\A^p$ \ is endowed with its \ $\theta-$structure. Writing the same exact sequence with the ordinary left-module structure of \ $\A^p$ \ gives
   \begin{equation*}
  0 \to \A^p \  \overset{^t\square.(Id_p.a - ^t\check{M})}{\longrightarrow} \ \A^p \to Ext^1_{\A}(E, \A) \to 0 . \tag{@@ bis}
  \end{equation*}
  where \ $^t\check{M}(b) : = \  ^tM(-b) $.\\
  Denote by \ $e^* : = (e_1^*, \cdots, e_p^*)$ \ the dual basis of \ $Hom_{\mathbb{C}[[b]]}(E,E_0)$. By definition of the action of \ $a$ \ on \ ${\it Hom}_{a,b}(E, E_0) $ \ we get, if \ $\omega$ \ is the class of 1 in \ $E_0$ :
  $$ (a.e_i^*)(e_j) = e_i^*(a.e_j) - a.e_i^*(e_j) = e_i^*(\sum_{h=1}^p \ m_{j,h}.e_h) - a.\delta_{i.j}.\omega = \check{m}_{j,i}.\omega $$ 
  because \ $a.\omega = 0$ \ in \ $E_0$, and the definition of the action of \ $b$ \ on \ ${\it Hom}_{a,b}(E, E_0) $. So we have \ $ a.e^* = ^t\check{M}.e^*$ \ concluding the proof. $\hfill \blacksquare$
  
  \begin{defn}\label{Dual} For any (a,b)-module \ $E$ \ the {\bf dual of \ $E$}, denoted by \ $E^*$, is the (a,b)-module \ $Ext^1_{\A}(E,\A) \simeq {\it Hom}_{a,b}(E, E_0) $.
  \end{defn}
  
  Of course, for any \ $\A-$linear map \ $f : E \to F$ \ between two (a,b)-modules we have an\ $\A-$linear ''dual'' map \ $f^* : F^* \to E^* $.\\
  It is an easy consequence of our previous description of \ $Ext^1_{\A}(E,\A)$ \ that we have a functorial isomorphism \ $(E^*)^*\to E$.
 
    \parag{Examples}
  \begin{enumerate}
  \item For each \ $\lambda \in \mathbb{C}$ \ we have \ $(E_{\lambda})^* \simeq E_{-\lambda}$.
  \item For \ $(\lambda,\mu) \in \mathbb{C}^2$ \ we have \ $E_{\lambda,\mu}^* \simeq E_{-\mu+1,-\lambda+1}$.
  \item Let \ $E$ \ be  the rank two simple pole (a,b)-module \ $E_1(0)$ \ defined by \ $a.e_1 = b.e_1 + b.e_2$ \ and \ $a.e_2 = b.e_2$.  Then its dual is isomorphic to \ $E_{-1}(0)$.\\
     It is also an elementary exercice to show the following isomorphisms :
    $$E_1(0) \simeq \mathbb{C}[[z]] \oplus  \mathbb{C}[[z]].Log z \quad {\rm and} \quad E_{-1}(0) \simeq  \mathbb{C}[[z]]\frac{1}{z^2} \oplus  \mathbb{C}[[z]].\frac{Log z}{z^2}$$ 
     with \ $a : = \times z$ \ and \ $ b : = \int_0^z $.
  \end{enumerate}
  
  \begin{prop}\label{Dualite et regularite}
  For any exact sequence of (a,b)-modules
  $$ 0 \to E' \overset{u}{\to} E \overset{v}{\to} E'' \to 0 $$
  we have an exact sequence of (a,b)-modules
  $$ 0 \to (E'')^* \overset{v^*}{\to}  E^* \overset{u^*}{\to}  (E')^* \to 0 .$$
  If \ $E$ \ is a simple p\^ole (a,b)-module, $E^*$ \ has a simple pole. \\
  For any regular (a,b)-module \ $E$ \ its dual \ $E^*$ \ is regular. Moreover, if \ $E^b$ \ and \ $E^{\sharp}$ \ are respectively the biggest simple pole submodule of \ $E$ \ and the saturation of \ $E$ \ by \ $b^{-1}.a$ \ in \ $E[b^{-1}]$, we have
  $$ (E^{\sharp})^* \simeq (E^*)^b \quad {\rm and} \quad  (E^b)^* \simeq (E^*)^{\sharp}.$$
  \end{prop} 
   
  \parag{Proof}
  The first assertion is a direct consequence of the vanishing of \ $Ext^i_{\A}(E, \A)$ \ for \ $i = 0, 2$, for any (a,b)-module and the long exact sequence for the ''Ext''.\\
  The condition that \ $E$ \ has a simple pole is equivalent to the fact that for any choosen basis \ $e$ \ of \ $E$ \ the matrix \ $M$ \ has its coefficients in \ $b.\A = \A.b$. Then this remains true for \ $^t\check{M}$.\\
  To prove the regularity of  \ $E^*$ \ when \ $E$ \ is  regular, we shall use induction on the rank of \ $E$. The rank 1 case is obvious because we have a simple pole in this case. Assume that the assertion is true for \ $ rank < p $ \ and consider a \ $rank = p$ \   regular (a,b)-module \ $E$. Using the theorem \ref{J-H} we have an exact sequence of (a,b)-modules
  $$ 0 \to E_{\lambda} \to E \to F \to 0  $$
  where \ $F$ \ is regular of rank \ $p-1$. This gives a short exact sequence
      $$ 0 \to F^* \to E^* \to E_{-\lambda} \to 0 $$
  and the regularity of \ $F^*$ \ and of \ $E_{-\lambda}$ \ implies the regularity of \ $E^*$.\\
  Now the inclusions \ $ E^b \subset E \subset E^{\sharp}$ \ gives exact sequences
  \begin{align*}
  & 0 \to Ext^1_{\A}(E/E^b, \A) \to E^* \to (E^b)^* \to Ext^2_{\A}(E/E^b, \A) \to 0 \\
  &  0 \to Ext^1_{\A}(E^{\sharp}/E, \A)  \to (E^{\sharp})^* \to E^* \to  Ext^2_{\A}(E^{\sharp}/E, \A) \to 0
  \end{align*}
  and the next lemma will show that the \ $Ext^1_{\A}(V, \A) = 0 $ \ for any \ $\A-$module which is a finite dimensional vector space, and also the finiteness (as a vector space) of \ $Ext^2_{\A}(V, \A)$.
  This implies  that we have, for any regular (a,b)-module, the inclusions
  $$ E^* \subset  (E^b)^*\quad {\rm and} \quad (E^{\sharp})^* \subset E^* .$$
  They imply, thanks to the fact that \ $(E^b)^*$ \ and \ $(E^{\sharp})^*$ \ have simple poles,
  $$  (E^*)^{\sharp} \subset (E^b)^* \quad {\rm and} \quad (E^{\sharp})^* \subset (E^*)^b.$$
    But the inclusion \ $(E^*)^b \subset E^* $ \ gives
     $$ E = (E^*)^*\subset ((E^*)^b)^* \subset ((E^{\sharp})^*)^* = E^{\sharp}$$
     and the minimality of \ $E^{\sharp}$ \ gives \ $((E^*)^b)^* = E^{\sharp}$ \ because \ $((E^*)^b)^*$ \ has a simple pole and contains \ $E$. Dualizing again  gives \ $(E^{\sharp})^* \simeq (E^*)^b$. 
     The last equality is obtained in a similar way from \ $ E^* \subset (E^*)^{\sharp}$. $\hfill \blacksquare$

  \begin{lemma}\label{Finite dim. modules}
  Let \ $V$ \ be an \ $\A-$module of finite dimension over \ $\mathbb{C}$. Then we have \ $Ext^i_{\A}(V,\A) = 0 $ \ for \ $i = 0,1$ \ and \ $Ext^2_{\A}(V,\A)$ \ is again a \ $\A-$module (via the \ $\theta-$structure of \ $\A$) which is a finite dimensional vector space. Moreover it has the same dimension than \ $V$ \ and there is a canonical \ $\A-$module isomorphism
  $$ Ext^2_{\A}(Ext^2_{\A}(V,\A),\A) \simeq V .$$
  \end{lemma}
  
  \parag{proof} We begin by proving the first assertion of the lemma for the special case \ $ V_{\lambda} : = \A \big/ \A.(a-\lambda) + \A.b$ \ for any  \ $\lambda \in \mathbb{C}$. Let us show that we have the free resolution
  $$ 0 \to \A \overset{\alpha}{\to} \A^2 \overset{\beta}{\to} \A \to V_{\lambda} \to 0 $$
  where \ $\alpha(x) : = (x.b, -x.(a-b-\lambda)), \quad \beta(u,v) : = u.(a-\lambda) + v.b $. The map \ $\alpha$ \ is clearly injective and \ $\beta(\alpha(x)) = x.(b.a - \lambda.b - (a-b-\lambda).b) = 0$. If we have \ $\beta(u,v) = 0$ \ then \ $u \in \A.b $; let \ $u = x.b$. Then we get
  $$ x.(a-b-\lambda).b + v.b = 0 \quad {\rm and \ so} \quad v = -x.(a-b-\lambda) .$$
  This gives the exactness of our resolution.\\
  Now the \ $Ext^i_{\A}(V_{\lambda},\A)$ \ are given by the cohomology of the complex
  $$ 0 \to \A   \overset{\beta^*}{\to}  \A^2  \overset{\alpha^*}{\to} \A \to 0 .$$
  The map \ $\beta^*(x) = ((a-\lambda).x, b.x)$ \ and \ $\alpha^*(u,v) = b.u - (a-b-\lambda).v $ \ are \ $\A-$linear for the \ $\theta-$structure of \ $\A$. Clearly \ $\beta^*$ \ is injective and \ $\alpha^*( \beta^*(x)) \equiv 0$. If \ $\alpha^*(u,v) = 0 $ \ set \ $v = b.y$ \ and conclude that \ $u = (a-\lambda).y$. This gives the vanishing of the \ $Ext^i$ \ for \ $i = 0,1$. The \ $Ext^2$ \ is the cokernel of \ $\beta^*$ \ which is easily seen to be isomorphic to \ $V_{\lambda}$.\\
  Consider now any finite dimensional \ $\A-$module \ $V$ \ over \ $\mathbb{C}$. We make an induction on \ $\dim_{\mathbb{C}}(V)$ \ to prove the vanishing of the \ $Ext^i$ \ for \ $i = 0,1$ \ and the assertion on the dimension of the \ $Ext^2$.\\
   The \ $\dim V = 1 $ \ case is clear because reduced to the case \ $V = V_{\lambda}$ \ for some \ $\lambda \in \mathbb{C}$. Assume that the case \ $\dim V = p$ \ is proved, for \ $p \geq 1$ \ and consider some \ $V$ \ with \ $\dim V = p +1$. Then \ $Ker \, b$ \ is not \ $\{0\}$ \ and is stable by \ $a$. Let \ $\lambda \in \mathbb{C}$ \ an eigenvalue of \ $a$ \ acting on \ $ Ker \, b$. Then a eigenvector generates in \ $V$ \ a sub-$\A-$module isomorphic to \ $V_{\lambda}$. \\
  The exact sequence of \ $\A-$modules
  $$ 0 \to V_{\lambda} \to V \to W \to 0 $$
  where \ $W : = V\big/ V_{\lambda}$ \ has dimension \ $p$ \ allows us to conclude, looking at the long exact sequence of  Ext .\\
  The last assertion follows from the remark that we produce a free resolution of \ $Ext^2_{\A}(V,\A)$ \ by taking \ $Hom_{\A}(-,\A)$ \ of a free (length two, see [B.97]) resolution of \ $V$ \ because of the already proved vanishing of the \ $Ext^i$ \  for \ $i = 0,1$. Doing this again gives back the initial resolution (remark that we use here that the \ $\theta\circ\theta-$structure on \ $Hom_{\A}(Hom_{\A}(\A,\A),\A)$ \ is the usual left structure on \ $\A$). $\hfill \blacksquare$
  
  \begin{cor}\label{Sym. Spec.}
 For a simple pole (a,b) module  \ $E$ \ denote by \ $S(E)$ \ the spectrum of \ $b^{-1}.a$ \ acting on \ $E/b.E$. Then we have
 $$ S(E^*) = - S(E).$$
 \end{cor}
 
 \parag{Proof} We make an induction on the rank of \ $E$. In rank $1$ \ the result is clear because we have \ $E \simeq E_{\lambda}$ \ for some \ $\lambda \in \mathbb{C}$, and \ $S(E_{\lambda}) = \{\lambda \}$. But we know that \ $E_{\lambda}^* = E_{-\lambda}$.\\
 Assume the assertion proved for any rank \ $p \geq 1$ \ simple pole (a,b)-module, and consider \ $E$ \ with rank \ $p +1$.  Using theorem \ref{J-H}, there exists \ $\lambda \in \mathbb{C}$ \ and an exact sequence  (a,b)-modules
 $$ 0 \to E_{\lambda} \to E \to F \to 0 $$
 where \ $rank(F)= p$ \ and where \ $F$ \ has a simple pole (because a quotient of a simple pole (a,b)-module has a simple pole !). The exact sequence of vector spaces 
 $$  0 \to E_{\lambda}/b.E_{\lambda}  \to E/b.E \to F/b.F \to 0 $$
 shows that \ $S(E) = S(F) \cup \{\lambda \}$. Now proposition \ref{Dualite et regularite} gives the exact sequence
 $$ 0 \to F^* \to E^* \to  E_{-\lambda} \to 0 $$
 which implies, as before, $S(E^*) = S(F^*) \cup \{-\lambda\} $. The induction hypothesis \ $S(F^*) = - S(F)$ \ allows to conclude. $\hfill \blacksquare$
 
 \bigskip
  
  \begin{lemma}\label{Duality for extension}
  For any pair of (a,b)-modules\ $E$ \ and\ $F$ \ there is a canonical isomorphism of vector spaces
  $$ D : Ext^1_{\A}(E,F) \to Ext^1_{\A}(F^*, E^*)  $$
 associated to the correspondance between 1-extensions (i.e. short exact sequences)
 $$ (0 \to F \to G \to E \to 0) \overset{D}{\to} (0 \to E^* \to G^* \to F^* \to 0).$$
 \end{lemma}
 
 \parag{Proof} We have a obvious isomorphism of \ $\mathbb{C}[[b]]-$modules\footnote{but be carefull with the \ $ b \to \check{b}$ !}
 $$ I :  Hom_b(E,F) \to Hom_b(Hom_b(F,E_0), Hom_b(E,E_0))\simeq Hom_b(F^*,E^*) $$
 because  \ $E_0 \simeq \mathbb{C}[[b]]$ \ as  a\ $\mathbb{C}[[b]]-$module. But recall that \ $Ext^1_{\A}(E,F)$ \ (resp. \ $Ext^1_{\A}(F^*,E^*)$) \ is the cokernel of the \ $\mathbb{C}-$linear map ''$a$'' defined on \ $ Hom_b(E,F) $ \ by the formula
 $$ (a.\varphi )(x) = \varphi(a.x) - a.\varphi(x) $$
 So it is enough to check that the isomorphism \ $I$ \ commutes with ''$a$'' in order to get an isomorphism between the cokernels of ''$a$'' in these two spaces.\\
 Let \ $\varphi \in Hom_b(E,F)$ and \ $\xi \in F^*$. Then \ $I(\varphi)(\xi) = \varphi\circ \xi $. 
 So, for \ $x \in  E$ \ we have (using \ $\Lambda$ \ to avoid too many ''$a$'')
 \begin{align*}
 & \Lambda(I(\varphi)(\xi) = I(\varphi)(a.\xi) - a.(I(\varphi)(\xi)) \\
 &  \Lambda(I(\varphi)(\xi)(x) = (\varphi \circ \xi )(a.x) - a.\xi(\varphi(x)) - \big(\xi(\varphi(a.x)) - a.\xi(\varphi(x))\big) \\
 & \quad\quad \quad\quad    = \big[(\Lambda(\varphi))\circ \xi \big](x) = I(\Lambda(\varphi))(x).
 \end{align*}
 So \ $\Lambda\circ I = I \circ \Lambda$. The map \ $I$ \ gives an isomorphism of complexes\\
 $$\xymatrix{Hom_{a,b}(E,F) \ar[r]^{\Lambda} \ar[d]^I & Hom_{a,b}(E,F)  \ar[d]^I \\
  Hom_{a,b}(F^*,E^*) \ar[r]^{\Lambda} &Hom_{a,b}(F^*,E^*) }$$
  and this conclude the proof, using lemma \ref{hom}. $\hfill \blacksquare$
  
  \bigskip
  
  For an (a,b)-module \ $E$ \ and an integer \ $m \in \mathbb{N}$ \ it is clear that \ $b^m.E$ \ is again an (a,b)-module. This can be generalize for any \ $m \in \mathbb{C}$.
  
  \begin{defn}
  For any (a,b)-module \ $E$ \ and any complex number \ $m \in \mathbb{C}$ \ define the (a,b)-module \ $b^m.E$ \ as follows :  as an \ $\mathbb{C}[[b]]-$module we let \ $b^m.E \simeq E \simeq \mathbb{C}[[b]]^{rank(E)}$;  the operator \ $a$ \ is defined as \ $a + m.b$.
  \end{defn}

  Precisely, this means that if \ $(e_1, \cdots, e_k)$ \ is a  \ $\mathbb{C}[[b]]-$basis of \ $E$ \ such that we have \ $a.e = M(b).e $ \ where \ $M \in End(\mathbb{C}^p) \otimes_{\mathbb{C}} \mathbb{C}[[b]]$, the (a,b)-module \ $b^m.E$ \ admit a basis, denote by \ $(b^m.e_1, \cdots, b^m.e_k)$, such that the operator \ $a$ \ is defined by the relation \ $a.(b^m.e) : = (M(b) + m.b.Id_k).(b^m.e) $.\\
  Remark that for \ $m \in \mathbb{N}$ \ this notation is compatible with the preexisting one, because of the relation \ $a.b^m = b^m.(a + m.b)$.\\
  For any \ $m \in \mathbb{N}$ \ there exists a canonical (a,b)-morphism
  $$ b^m.E \to  E  $$
  which is an isomorphism of \ $b^m.E$ \ on \ $Im(b^m : E \to E)$. But remark that the map \ $b^m : E \to E$ \ is not $a-$linear (but the image is stable by \ $a$).\\
  For any \ $m \in \mathbb{N}$ \ there is also a canonical (a,b)-morphism
  $$ E \to b^{-m}.E $$
  which induces an isomorphism of \ $E$ \ on \ $Im(b^m : b^{-m}.E \to b^{-m}.E)$. So we may write, via this canonical identification,  $b^m.(b^{-m}.E) = E$.\\
  It is easy to see that for any \ $m, m' \in \mathbb{C}$ \ we have a natural isomorphism
  $$ b^{m'}.(b^m.E) \simeq b^{m+m'}.E \quad {\rm and \ also} \quad b^0.E \simeq E .$$
  
  \parag{Remark} It is easy to show that for any \ $m \in \mathbb{C}$ \ there exists an unique $\mathbb{C}-$algebra automorphism 
   $$ \eta_m : \A \to \A \quad {\rm such \  that } \quad \eta(1) = 1, \eta(b) = b \quad {\rm and} \quad \eta(a) = a +m.b .$$
   Using this automorphism, one can define a left \ $\A-$module  \ $b^m.F$ \ for any left \ $\A-$module \ $F$ \ and any \ $m \in \mathbb{C}$. This is, of course compatible with our definition in the context of (a,b)-modules. $\hfill \square$
   
   \bigskip
   
   The behaviour of the correspondance \ $E \to b^m.E$ \  by duality is given by the following easy lemma; the proof is left as an exercice.
   
   \begin{lemma}
   For any (a,b)-module \ $E$ \ and any \ $m \in \mathbb{C}$ \ there is natural (a,b)-isomorphism
   $$ (b^m.E)^*\to b^{-m}.E^* .$$
   \end{lemma}
   
   \bigskip
   
   The following corollary of the lemma  \ref{inclusion} and the proposition \ref{Dualite et regularite} allows to show that duality preserves the index.
   
   \begin{lemma}\label{or dual}
   Let \ $E$ \ be a regular (a,b)-module. Then  we have \ $\delta(E^*) = \delta(E)$.
   \end{lemma}
   
   \parag{Proof} By definition \ $\delta(E)$ \ is the smallest integer \ $k \in \mathbb{N}$ \  such that \ $E^{\sharp} \subset b^{-k}.E$.\\
 Now \ $E^{\sharp} \subset b^{-m}.E$ \  implies by duality that \ $b^m.E^* \subset (E^*)^b$. So,  by lemma \ref{inclusion}, we have \ $m \geq \delta(E^*)$. This proves that \ $\delta(E) \leq \delta(E^*)$ \ and we obtain the equality by symetry. $\hfill \blacksquare$
 
\parag{Remark} Duality does not preserve the order of regularity : in the example given before the definition \ref{$E^b$} we have \  $or(E) = 2$ \ and \ $or(E^*) = 1$. $\hfill \square$

\bigskip
   
   Let us conclude this section by an easy exercice.
   
   \parag{Exercice}
   For any (a,b)-modules  \ $E, F$ \ and any \ $\lambda \in \mathbb{C}$ \ there are natural (a,b)-isomorphisms
   \begin{enumerate}
   \item  \ $b^{\lambda}.E_{\mu} \simeq E_{\lambda+\mu}$.
    \item \ $b^{\lambda}.Hom_{a,b}(E, F) \simeq  Hom_{a,b}(b^{-\lambda}.E, F) \simeq Hom_{a,b}(E, b^{\lambda}.F)$. 
    \item Then deduce from the previous isomorphisms that \ $Hom_{a,b}(E, E_{\lambda}) \simeq b^{-\lambda}.E^*$, and \ $Ext^1_{\A}(E, E_{\lambda}) \simeq E^*/(a + \lambda.b).E^*$.
    \end{enumerate}

\subsection{Width of a regular (a,b)-module.}

\parag{Notation} For a complex number \ $\lambda$ \ we shall note by \ $\tilde{\lambda}$ \ is class in \ $\mathbb{C}\big/\mathbb{Z}$. We shall order elements in each class \ $modulo \  \mathbb{Z}$ \ by its natural order on real parts. $\hfill \square$

\begin{defn}\label{largeur}
Let \ $E$ \ be a regular (a,b)-module and let \ $\tilde{\lambda} \in \mathbb{C}\big/\mathbb{Z}$. We define the following complex numbers :
\begin{align*}
& \tilde{\lambda}_{min}(E) : = \inf \{\lambda \in \tilde{\lambda} /\exists \ {\rm a \ non \ zero \ morphism} \quad E_{\lambda} \rightarrow E  \} \\
& \tilde{\lambda}_{max}(E) = \sup \{ \lambda \in \tilde{\lambda} / \exists  \ {\rm a \ non \ zero \ morphism} \quad  E \to E_{\lambda}  \} \\
& L_{\tilde{\lambda}}(E) = \tilde{\lambda}_{max}(E) - \tilde{\lambda}_{min}(E) \in \mathbb{Z} \\
& L(E) = \sup \{ \tilde{\lambda} \in  \mathbb{C}/\mathbb{Z} \  /  \ L_{\tilde{\lambda}}(E) \}
\end{align*}
with the following conventions :
 \begin{align*}
 & \inf \{\emptyset \} = +\infty, \ \sup \{\emptyset \} = - \infty \quad {\rm and} \\
 & -\infty - \lambda = - \infty \quad \forall \lambda \in ]-\infty, +\infty] \\
 & +\infty - \lambda = + \infty \quad \forall \lambda \in [-\infty, +\infty[. 
  \end{align*}
We shall call \ $L(E)$ \  {\bf the width of\ $E$}.
\end{defn}

\parag{Remarks}
\begin{enumerate}
\item A non zero morphism \ $E_{\lambda} \to E $ \ is necessarily injective. Either its image is a normal submodule in \ $E$ \ or there exists an integer \ $k \geq 1$ \ and a morphism \ $E_{\lambda-k} \to E$ \ whose image is normal an contains the image of the previous one.
\item In a dual way, a non zero morphism \ $E \to E_{\lambda}$ \  has an image equal to \ $b^k.E_{\lambda} \simeq E_{\lambda+k}$, where \ $k \in \mathbb{N}$.
\item A non zero morphism \ $E_{\lambda} \to E_{\mu}$ \ implies that \ $\lambda $ \ lies in \ $\mu + \mathbb{N}$. It is possible that for some \ $E$ \ we have \ $\tilde{\lambda}_{max}(E) < \tilde{\lambda}_{min}(E) $. For instance this is the case for the rank 2 regular (a,b)-module \ $E_{\lambda,\mu}$ \ from \ref{class}. So the width of a regular but not simple pole (a,b)-module  is not necessarily a non negative integer. 
\item Let \ $E$ \ and \ $F$ \ be regular (a,b)-modules. If there is a surjective morphism \ $ E \to F$ \ then for all \ $\tilde{\lambda} \in \mathbb{C}\big/\mathbb{Z}$ \ we have \ $\tilde{\lambda}_{max}(E) \geq \tilde{\lambda}_{max}(F)$. \\
If there is an injective morphism \ $E' \to E$ \ then for all \ $\tilde{\lambda} \in \mathbb{C}\big/\mathbb{Z}$ \ we have \ $\tilde{\lambda}_{min}(E) \leq \tilde{\lambda}_{min}(E')$.

\item Every submodule of \ $E$ \ isomorphic to \ $E_{\lambda}$ \ is contained in \ $E^b$. So we have \ $\tilde{\lambda}_{min}(E) = \tilde{\lambda}_{min}(E^b)$, for every regular (a,b)-module \ $E$ \ and every  \ $\tilde{\lambda} \in \mathbb{C}\big/\mathbb{Z}$.
\item In a dual way, every morphism \ $E \to E_{\lambda}$ \ extends uniquely to a morphism \ $E^{\sharp} \to  E_{\lambda}$ \ with the same image. So for every regular (a,b)-module \ $E$ \ and every  \ $\tilde{\lambda} \in \mathbb{C}\big/\mathbb{Z}$, we get \ $\tilde{\lambda}_{max}(E) = \tilde{\lambda}_{max}(E^{\sharp})$. $\hfill \square$
\end{enumerate}

\begin{lemma}\label{obvious}
\begin{enumerate}
\item Let \ $E$ \  a simple pole (a,b)-module and let \ $S(E)$ \ denotes the spectrum of the linear map \  $b^{-1}.a : E/b.E \to E/b.E$, we have
\begin{equation*}
 \tilde{\lambda}_{min}(E) = \inf \{ \lambda \in S(E) \cap \tilde{\lambda} \} \quad {\rm and} \quad   \tilde{\lambda}_{max}(E) = \sup \{ \lambda \in S(E) \cap \tilde{\lambda} \} \tag{@}
 \end{equation*}
\item For any regular (a,b)-module \ $E$ \ we have
$$ \widetilde{(-\lambda)}_{max}(E^*) = - \tilde{\lambda}_{min}(E) \quad\qquad \widetilde{(-\lambda)}_{min}(E^*) = -  \tilde{\lambda}_{max}(E).$$
This implies \ $ L_{-\tilde{\lambda}}(E^*) = L_{\tilde{\lambda}}(E)\  \forall \tilde{\lambda} \in \mathbb{C}/\mathbb{Z}$, and so \ $L(E^*) = L(E)$.
\item For any regular (a,b)-module \ $E$ \ and any \ $\tilde{\lambda} \in \mathbb{C}\big/\mathbb{Z}$ \ we have equivalence between
$$\tilde{\lambda}_{min}(E)  \not= +\infty \quad {\rm and} \quad \tilde{\lambda}_{max}(E) \not= -\infty .$$
\end{enumerate}
\end{lemma}

\parag{Proof} Let \ $E$ \ be a simple pole (a,b)-module. We have already seen (in proposition \ref{sub min})  that if \ $\lambda \in S(E)$ \ is minimal in its class modulo $1$, there exists a non zero \ $x \in E$ \ such that \ $a.x = \lambda.b.x$. This implies that \ $\tilde{\lambda}_{min} \leq \inf \{ \lambda \in S(E) \cap \tilde{\lambda} \}$. But the opposite inequality is obvious, so the first part of  (@) is proved.\\
Using corollary \ref{Sym. Spec.} and the result already obtained  for \ $E^*$ \ gives
$$ \widetilde{(-\lambda)}_{min}(E^*) = \inf \{-\lambda \in S(E^*) \cap \widetilde{(-\lambda)} \} = - \sup \{\lambda \in S(E) \cap \tilde{\lambda}\}.$$
So for \ $\mu = \sup \{\lambda \in S(E) \cap \tilde{\lambda}\}$ \ we have an exact sequence of (a,b)-modules
$$ 0 \to E_{-\mu} \to E^* \to F \to 0 $$
and by duality, a surjective map \ $E \to E_{\mu}$. This implies \ $\tilde{\lambda}_{max} \geq \mu$. As, again, the opposite inequality is obvious, the second part of (@) is proved.\\
Let us prove now the relations in 2. \\
Remark first that these equalities are true for a simple pole (a,b)-module because of \ $(@)$ \ and corollary \ref{Sym. Spec.}.\\
For any regular (a,b)-module \ $E$ \ we know that 
 $$\tilde{\lambda}_{min}(E) =  \tilde{\lambda}_{min}(E^b) =  \inf \{ \lambda \in S(E^b) \cap \tilde{\lambda} \} \quad {\rm and} \quad  \widetilde{(-\lambda)}_{max}(E^*) = \widetilde{(-\lambda)}_{max}((E^*)^{\sharp}).$$
  But we have
   $$  \widetilde{(-\lambda)}_{max}((E^*)^{\sharp}) = \sup\{-\lambda \in S((E^*)^{\sharp}) \cap \widetilde{(-\lambda)} \} = - \inf\{\lambda \in S((E^*)^{\sharp})^*\cap \tilde{\lambda} \} $$ 
   because \ $(E^*)^{\sharp}$ \  has a simple pole, using corollary \ref{Sym. Spec.}. So we obtain
   $$ \widetilde{(-\lambda)}_{max}(E^*) = - \tilde{\lambda}_{min}(E^b) = - \tilde{\lambda}_{min}(E)$$
   because \ $(E^*)^{\sharp})^* = E^b$ \ (see proposition \ref{Dualite et regularite}).\\  
The second relation is analoguous.\\
The equivalence in 3  is obvious in the simple pole case using \ $(@)$.\\
 The general case is an easy consequence using \ $E^b, E^{\sharp}$ : if \ $\tilde{\lambda}_{min}(E) \not= +\infty$ \ so is \ $\tilde{\lambda}_{min}(E^{\sharp})$ \ because \ $E \subset E^{\sharp}$. Then \ $\tilde{\lambda}_{max}(E^{\sharp})\not= -\infty$ \ and so is \ $\tilde{\lambda}_{max}(E)$. The converse is analoguous using \ $E^b$. \hfill $\blacksquare$

\parag{Remarks}
\begin{enumerate}
\item If \ $E$ \ has a simple pole, we have   \ $L_{\tilde{\lambda}}(E) \geq 0$ \ or \ $L_{\tilde{\lambda}}(E) = - \infty$ \ for any \ $\tilde{\lambda}$ \ in \ $\mathbb{C}/\mathbb{Z}$. So \ $L(E)$ \ is always \ $\geq 0$.
\item In cases 1 and 2 of the proposition \ref{class} the formula \ $(@)$ \ gives the values of \ $\tilde{\lambda}_{min}$ \ and \ $\tilde{\lambda}_{max}$ \ for any \ $\tilde{\lambda} \in \mathbb{C}/\mathbb{Z}$.\\
For the remaining cases we can compute these numbers using the fact that we already know the corresponding \ $E^b$ \ and \ $E^{\sharp}$ \ and the remark 5 and 6  before the preceeding lemma. $\hfill \square$
\end{enumerate}

\bigskip

\begin{prop}\label{Induction largeur}
Let \ $E$ \ be a regular (a,b)-module and let \ $\tilde{\lambda} \in \mathbb{C}\big/\mathbb{Z}$. Assume that \ $ \lambda = \tilde{\lambda}_{min}(E) < +\infty $. Consider an exact sequence of (a,b)-modules
$$ 0 \to E_{\lambda} \to E \overset{\pi}{\to} F \to 0 .$$
Then we have for all  \ $\tilde{\mu} \in \mathbb{C}/\mathbb{Z}$ \ the inequality
\begin{equation*}
 L_{\tilde{\mu}}(F) \leq  L_{\tilde{\mu}}(E) + 1. \tag{i}
 \end{equation*}
 \end{prop}
 
 \parag{Proof} As \ $\tilde{\mu}_{max}(F) \leq \tilde{\mu}_{max}(E)$ \ for any \ $\mu \in \mathbb{C}$ \ it is enough to prove that we have \ $ \tilde{\mu}_{min}(E) \leq \tilde{\mu}_{min}(F) +1$ \ for all \ $\tilde{\mu} \in \mathbb{C}/\mathbb{Z}$.
 
  Let begin by the case of \ $\tilde{\mu} = \tilde{\lambda}$. We want to show the inequality
 \begin{equation*}
 \tilde{\lambda}_{min}(F) \geq \lambda -1 \tag{ii}
 \end{equation*}
 Let  \ $E_{\lambda - d} \hookrightarrow F$ \ with \ $d \geq 0$. The rank 2 (a,b)-module  \ $G : = \pi^{-1}(E_{\lambda - d})$ \ is contained in \ $E$, so \ $\lambda = \tilde{\lambda}_{min}(G)$. We have
  the exact sequence of (a,b)-modules
 $$ 0 \to E_{\lambda} \to \pi^{-1}(E_{\lambda - d}) \overset{\pi}{\to} E_{\lambda-d} \to 0. $$
 Now let us  compare \ $G$ \ with the list in proposition \ref{class}.\\
  If \ $G$ \ is in case 1, we have \ $E_{\lambda-d} \subset G$ \ so \ $d = 0$ \ because \ $\lambda = \tilde{\lambda}_{min}(G)$.\\
 If \ $G$ \ is in case 2, we have \ $\lambda - d = \lambda + n$ \ with \ $n \in \mathbb{N}$, so \ $d = 0$.\\
 If \ $G$ \ is in case 3, we have  \ $G \simeq E_{\lambda, \lambda+k}$ \ with \ $k \in \mathbb{N}$. Then  the theorem \ref{J-H} gives \ $2\lambda - d = 2\lambda +k -1 $ \ and so \ $d = 1-k \leq 1$.\\
 If \ $G$ \ is in case 4, we have \ $G \simeq E_{\lambda,\lambda+n}(\alpha)$. Again theorem \ref{J-H} gives \ $2\lambda - d = 2\lambda+n-1 $ \ so \ $d = 1-n \leq 0$ \ because \ $n \in \mathbb{N}^*$. So \ $d = 0$.\\
  We conclude that we always have \ $d \leq 1$ \ and this proves (ii).\\
   
    \smallskip
   
For \ $\tilde{\mu} \not= \tilde{\lambda}$ \ let us prove now  the following inequality :
\begin{equation*}
\tilde{\mu}_{min}(F) \leq  \tilde{\mu}_{min}(E) \leq \tilde{\mu}_{min}(F) +1 . \tag{iii}
\end{equation*}
Consider an injective morphism  \ $E_{\mu} \to E$ \ with \ $\mu = \tilde{\mu}_{min}(E)$. The restriction of \ $\pi$ \ to \ $E_{\mu}$ \ is injective and so it gives \  $ \tilde{\mu}_{min}(E) \geq \tilde{\mu}_{min}(F)$. Assume now that we have an injective morphism  \ $E_{\mu'} \hookrightarrow F$ \ with  \ $\mu' = \tilde{\mu}_{min}(F)$, and consider the rank 2 (a,b)-module \ $\pi^{-1}(E_{\mu'})$. Using the proposition \ref{class} where only cases 1 or 3 are possible now, it can be easily check that (iii) is satisfied. $\hfill \blacksquare$

\parag{Remarks} 
\begin{enumerate}
\item In the situation of the previous proposition we have either \ $\tilde{\lambda}_{min}(E) \geq \tilde{\lambda}_{max}(E)$ \ or \ $\tilde{\lambda}_{max}(E)= \tilde{\lambda}_{max}(F)$ : Assume that we have  \ $\lambda < \lambda' : = \tilde{\lambda}_{max}(E)$. Then there exists a surjective morphism \ $ q : E \to E_{\lambda'}$, and, as the restriction of \ $q$ \ to  \ $E_{\lambda}$ \ is zero, the map \ $q$ \ can be factorized and gives a surjective morphism \ $\tilde{q} : F \to E_{\lambda'}$. So we get \ $\tilde{\lambda}_{max}(E)\leq \tilde{\lambda}_{max}(F)$, and the desired equality thanks to the preceeding lemma.
\item We shall use later that in the situation of the previous proposition we have the inequality \ 
$ \tilde{\lambda}_{max}(F) \leq \lambda + L(E) . \hfill  \square$
\end{enumerate}

\begin{cor}\label{maj. lambda(max)(F)}
In the situation of the previous proposition we have the inequality  \ $L(E) + rank(E) \geq L(F) + rank(F)$. So this integer is always positive for any non zero regular (a,b)-module.
\end{cor}

\parag{Proof} As  the rank  1 case is obious, an easy induction on the rank of \ $E$ \ using the propositions \ref{sub min}  and \ref{Induction largeur} gives the proof. $\hfill \blacksquare$

\parag{Examples}
\begin{enumerate}
\item The (a,b)-module
$$  J_k(\lambda) : = \A\big/\A.(a -(\lambda+k-1).b)(a -(\lambda+k-2).b)\cdots (a-\lambda.b)$$
which has rank \ $k$, satisfies \ $\lambda_{max} = \lambda$ \ and \ $\lambda_{min} = \lambda+k-1$. So its width is \ $L(J_k(\lambda)) = -k + 1 $.\\
To understand easily the (a,b)-module \ $J_k(\lambda)$ \ the reader may use the following alternative definition of it : there is a \ $\mathbb{C}[[b]]-$basis  \ $(e_1, \cdots, e_k)$ \ in which the action of ''$a$'' is given by 
$$ a.e_1 = e_2 + \lambda.b.e_1,\quad  a.e_2 = e_3 + (\lambda+1).b.e_2, \cdots, a.e_k =(\lambda+ k-1).b.e_k .$$
\item The rank 2 (a,b)-module \ $E_{\lambda} \oplus E_{\lambda+n}$ \ has width \ $n$. This shows that, despite the fact that the width is always bigger than \ $- rank(E) + 1$, the width may be arbitrarily big, even for a rank 2  regular (a,b)-module. $\hfill \square$
\end{enumerate}

 \section{Finite determination of regular (a,b)-modules.}

\subsection{Some more preliminaries.}

 \begin{lemma} \label{finite det. reg.}
 Let \ $E$ \ be a regular (a,b)-module of index \ $\delta(E) = k$. For \ $N \geq k+1$ \ the quotient map \ $ q_N : E \to E\big/b^N.E $ \ induces a bijection between simple pole sub-(a,b)-modules \ $F$ \ containing \ $b^k.E^{\sharp}$ \ and sub \ $\A-$modules \ $\mathcal{F} \subset  E\big/b^N.E $ \ satisfying the following two conditions
 \begin{enumerate}[i)]
 \item \ $ a.\mathcal{F} \subset b.\mathcal{F}$ ;
 \item  \ $b^k.E^{\sharp}\big/b^N.E \subset \mathcal{F}$.
 \end{enumerate}
 \end{lemma}
 
 \parag{Proof} It is clear that if \ $F$ \ is a simple pole sub-(a,b)-module of \ $E$ \ containing \ $b^k.E^{\sharp}$ \  the image \ $\mathcal{F} : = q_N(F)$ \ is a \ $\A-$submodule of \ $E\big/b^N.E $ \ such that i) and  ii) are fullfilled. Conversely, if a \ $\A-$submodule \ $\mathcal{F}$ \ satisfies i) and ii), let \ $F : = q_N^{-1}(\mathcal{F})$. Of course, \ $F$ \ is a sub-(a,b)-module of \ $E$ \ and contains \ $b^k.E^{\sharp}$. The only point to see is that \ $F$ \ has a simple pole. If \ $x \in F$ \ then \ $a.q_N(x) \in b.\mathcal{F}$ \ so \ $a.x \in b.F + b^N.E$. As \ $N \geq k+1$ \ we may write \ $ a.x = b.y + b.z $ \ with \ $y \in F $ \ and \ $z \in b^{N-1}.E \subset  b^k.E^{\sharp} \subset F $. This completes the proof. $\hfill \blacksquare$.
 
 \parag{Remarks} 
 \begin{enumerate}
 \item we may replace \ $b^k.E^{\sharp}$ \ by \ $b^k.E$ \ in the second condition imposed on \ $F$ \ and \ $\mathcal{F}$ : if a simple pole (a,b)-submodule \ $F$ \ contains \ $b^k.E$ \ it contains \ $b^k.E^{\sharp}$ \ by definition of \ $E^{\sharp}$. This allows to avoid the use of \ $E^{\sharp}$ \ in the previous lemma.
 \item The biggest \ $\mathcal{F}$ \ satisfying i) and ii) corresponds to \ $E^b$. So we may recover \ $E^b$ \ from the quotient \ $E\big/b^N.E$ \ for \ $N \geq \delta(E)+1$. $\hfill \square$
 \end{enumerate}
 
 \begin{cor}\label{det. finie ordre reg.}
  Let \ $E$ \ be a regular (a,b)-module of order of regularity \ $k$. Fix \ $N \geq k+1$ \ and assume that we has an isomorphism of \ $\A-$modules 
   $$ \varphi : E\big/b^N.E \to  E'\big/b^N.E'  $$
   where \ $E'$ \ is an (a,b)-module. Then \ $E'$ \ is regular and has order of regularity \ $k$. Moreover we  have the equality \ $\varphi(E^b\big/b^N.E) = (E')^b\big/b^N.E' $.
  \end{cor}
  
  \parag{Proof} As \ $k$ \ is the order of regularity of \ $E$ \ we have \ $a^{k+1}.E \subset \sum_{j=0}^{k} \   a^j.b^{k-j+1}E $. The inequality \ $N \geq k+1$ \ gives \ $a^{k+1}.E\big/b^N.E \subset \sum_{j=0}^{k} \  a^j.b^{k-j+1}E\big/b^N.E $, and this is also true for \ $E'\big/b^N.E'$, and implies \ $a^{k+1}.E' \subset \sum_{j=0}^{k} \  a^j.b^{k-j+1}E' $. So the order of regularity of \ $E'$ \ is at most \ $k$. We conclude that it is exactly \ $k$ \ by symetry.\\
  The last stament comes from the second remark above, as \ $or(E) \geq \delta(E)$. $\hfill \blacksquare$
  
\subsection{Finite determination for a rank one extension.}

 \begin{lemma}\label{existence N}
Let \ $E$ \ be an (a,b)-module et  fix a complex number \ $\lambda$. There exists \ $N(E,\lambda) \in \mathbb{N}$ \ such that for any \ $N \geq N(E,\lambda)$ we have the following inclusion :
  $$b^N.E \subset (a - \lambda.b).E .$$
\end{lemma}

\parag{Proof} With the $b-$adic topology, \ $E$ \ is a Frechet space. The \ $ \mathbb{C}-$linear map \ $a - \lambda.b : E \to E $ \ is continuous. The finiteness theorem of [B.95], theorem 1.bis p.31 gives that the kernel and cokernel of this map are finite dimensional vector spaces. So the subspace \ $(a - \lambda.b).E$ \ is closed in \ $E$. This statement corresponds to the equality
\begin{equation*}
 \cap_{N \geq 0} \big[(a - \lambda.b).E + b^N.E \big] = (a- \lambda.b).E  \tag{@}
 \end{equation*}
 But the images of the subspaces \ $b^N.E$ \ in the finite dimensional vector space \\
  $E \big/ (a- \lambda.b).E$ \ is a decreasing sequence. So it is stationnary, and, as the intersection is \ $\{0\}$ \ thanks to \ $(@)$, the result follows. $\hfill \blacksquare$

  \begin{prop}\label{4}
   Let \ $F$ \ be an (a,b)-module and \ $\lambda$ \ a complex number. Consider a short exact sequence of (a,b)-modules
  \begin{equation*}
   0 \to E_{\lambda} \overset{\alpha}{\to} E \overset{\beta}{\to} F \to 0  \tag{$@@$}
   \end{equation*}
  where \ $E_{\lambda} : = \A\big/\A.(a - \lambda.b) $. Then, for any \ $N \geq N(F^*, -\lambda)$, the extension \ $(@@)$ \ is uniquely determined by the following extension of \ $\A-$modules which are finite dimensional vectors spaces
  \begin{equation*}
   0 \to  E_{\lambda}\big/b^N.E_{\lambda}  \overset{\alpha}{\to} E\big/b^N.E \overset{\beta}{\to} F\big/b^N.F  \to 0  \tag{$@@_N$}   \end{equation*}
 obtained from \ $(@@)$ \ by ''quotient by \ $b^N$''.
\end{prop}

\parag{Comments} This statement needs some more explanations. Denote by \ $K_N$ \ the kernel of the obvious map (forget "a") 
 $$ob_N : Ext^1_{\A}(F/b^N.F, E_{\lambda}/b^N.E_{\lambda}) \to Ext^1_b(F/b^N.F, E_{\lambda}/b^N.E_{\lambda})$$
 where \ $Ext^1_b(-,-)$ \ is a short notation for \ $Ext^1_{\mathbb{C}[[b]]}(-,-)$. The  short exact sequence correspondance \ $(@@) \to (@@_N)$ \ gives a map 
 $$\delta_N : Ext^1_{\A}(F, E_{\lambda}) \to Ext^1_{\A}(F/b^N.F, E_{\lambda}/b^N.E_{\lambda}) $$
 whose  range  lies in \ $K_N$, because the \ $\mathbb{C}[[b]]-$exact sequence \ $(@@)$ \ is split as \ $F$ \ is \ $\mathbb{C}[[b]]-$free, and so is the exact sequence \ $(@@_N)$. The precise signification of the previous proposition is that for \ $N \geq N(F^*, -\lambda)$ \ the map \ $\delta_N$ \ is a \ $\mathbb{C}-$linear  isomorphism between the vector spaces  \ $Ext^1_{\A}(F, E_{\lambda})$ \ and \ $K_N$. $\hfill \square$
 
  \parag{Proof} As a first step to realize the map \ $\delta_N$ \ let us consider  the following commutative  diagramm of complex vector spaces, deduced from the exact sequences of \ $\A-$modules:
  \begin{align*}
 &  0 \to E_{\lambda + N} \to E_{\lambda} \to E_{\lambda}/b^N.E_{\lambda} \to 0  \\
  & 0 \to b^N.F \to F \to F/b^N.F \to 0 
  \end{align*}

 $$ \xymatrix{ Ext^1( F/ b^N.F,E_{\lambda+N}) \ar[d] \ar[r] &  Ext^1(F,E_{\lambda+N})  \ar[d]^{\alpha}  \ar[r] &  Ext^1(b^N.F,E_{\lambda +N})  \ar[d] \\
  Ext^1( F/ b^N.F, E_{\lambda}) \ar[d] \ar[r] & Ext^1(F, E_{\lambda})  \ar[d]^{\beta} \ar[r]^u  & Ext^1(b^N.F,E_{\lambda})  \ar[d]^v \\
  Ext^1(F/ b^N.F,E_{\lambda}/b^N.E_{\lambda}) \ar[r]^>>>>{\gamma} &  Ext^1(F, E_{\lambda}/b^N.E_{\lambda}) \ar[r]^w  & Ext^1(b^N.F, E_{\lambda}/b^N.E_{\lambda})} $$
  
  We have the following propreties : 
  
  \begin{enumerate}
 \item The surjectivity of the map \ $\beta$ \  is consequence of the vanishing of the vector space \ $Ext^2_{\A}(F, E_{\lambda +N})$ \ thanks to the proposition \ref{Dualite et regularite}.
  \item the vanishing of the composition \ $u\circ v $ \  is consequence of lemma  \ref{hom} and of the fact that the restriction map 
  $$ Hom_b(F, E_{\lambda}) \to Hom_b(b^N.F, E_{\lambda}) \to Hom_b(b^N.F, E_{\lambda}/b^N.E_{\lambda})$$ 
   is obviously zero.
   \item So the map \ $w$ \ is zero and \ $\gamma$ \ is surjective.
   \item The kernel of \ $\gamma$ \ is given by the image of the injective map
 $$ \partial : Hom_{\A}(b^N.F, E_{\lambda}/b^N.E_{\lambda}) \hookrightarrow Ext^1_{\A}(F/ b^N.F,E_{\lambda}/b^N.E_{\lambda}) .$$ 
 This  is  a consequence of the vanishing of the map
 $$ Ext^0_{\A}(F, E_{\lambda}/b^N.E_{\lambda}) \to Ext^0_{\A}(b^N.F, E_{\lambda}/b^N.E_{\lambda}).$$
 \end{enumerate}
 Let us show now that for \ $N \geq N(F^*, -\lambda)$ \ the map \ $\alpha$ \ is zero.  Using again the isomorphisms given by the lemma \ref{hom}, \ $\alpha$ \ is induced by the obvious map \ $ Hom_b(F, b^N.E_{\lambda}) \to Hom_b(F, E_{\lambda})$, whose image is \ $b^N.Hom_b(F, E_{\lambda})$. Denote respectively by \ $G$ \ and \ $H$ \  the (a,b)-modules given by \ $Hom_b(F, b^N.E_{\lambda})$ \ and \ $Hom_b(F, E_{\lambda})$ \ with the action of "$a$" defined by \ $\Lambda$ (see \ref{hom}). Then we have the following commutative diagramm
   $$ \xymatrix{ G \ar[r]^{i} \ar[d] & H \ar[d] \\ G/a.G \ar[r] \ar[d]^{\simeq} & H/a.H \ar[d]^{\simeq}\\
   Ext^1_{\A}(F, b^N.E_{\lambda}) \ar[r]^{\alpha} & Ext^1_{\A}(F, E_{\lambda})}$$
   and the image of \ $i$ \ is \ $b^N.H$. So the map \ $\alpha$ \ will be zero as soon as \ $b^N.H \subset a.H$ \ and this is fullfilled for \ $N \geq  N(H,0) = N(F^*, -\lambda) $. This last equality coming from the isomorphisms
   $$H/a.H \simeq Ext_{\A}^1(F, E_{\lambda}) \simeq Ext_{\A}^1( E_{-\lambda}, F^*) \simeq F^*/(a+\lambda.b).F^*$$
   see the exercice concluding section 3.3.
   
   \smallskip
   
  Consider now the commutative diagramm
   
   $$ \xymatrix{\quad & 0 \ar[d] & K_N \ar[d]^i &Ext_{\A}^1(F, E_{\lambda})  \ar[d]^{\beta} \ar[l]_{\hat{\delta}_N} \ar[ld]_{\delta_N}  \\
    0 \ar[r] & Hom_{\A}(b^N.F, E_{\lambda}/b^N.E_{\lambda}) \ar[r]^{\partial}\ar[d]^{ob_N} &  Ext^1_{\A}(F/ b^N.F,E_{\lambda}/b^N.E_{\lambda})\ar[r]^>>>>>{\gamma} \ar[d]^{ob_N}& Ext^1_{\A}(F, E_{\lambda}/b^N.E_{\lambda})\\ \quad & Hom_b(b^N.F, E_{\lambda}/b^N.E_{\lambda})\ar[r]^{\simeq} &  Ext^1_b(F/ b^N.F,E_{\lambda}/b^N.E_{\lambda})& \quad }$$
    The surjectivity of \ $\beta$ \ implies that the map \ $i\circ\gamma $ \ is surjective ( we know that the extensions in the image of \ $\delta_N$ \ comes from \ $K_N$, so \ $\delta_N$ \ factors in \ $\hat{\delta}_N\circ i$).\\
    We have \ $i(K_N) \cap Im(\partial_N) = (0)$ \ because \ $ob_N$ \ is injective on \ $Im(\partial_N)$. \\
    So \ $i$ \ induces an isomorphism of vector spaces from \ $K_N$ \ to 
     $$ Ext^1_{\A}(F/ b^N.F,E_{\lambda}/b^N.E_{\lambda})/ Im(\partial_N)\overset{\gamma}{ \simeq} Ext^1_{\A}(F, E_{\lambda}/b^N.E_{\lambda}) \overset{\beta^{-1}}{\simeq} Ext_{\A}^1(F, E_{\lambda}).$$
     This completes the proof . $\hfill \blacksquare$
     
     \bigskip
     We shall need some bound for the integer \ $N(F^*, -\lambda)$ \ which appears in the previous proposition for the proof of our theorem.
      
      \begin{lemma}
      Let \ $G$ \ be a regular (a,b)-module and let \ $\mu \in \mathbb{C}$. A sufficient condition on \ $N \in \mathbb{N}$ \ in order to have the inclusion \ $ b^N.G \subset (a-\mu.b).G$ \ is  
      $$N \geq \mu -\tilde{\mu}_{min}(G) + \delta(G) + 2.$$
      \end{lemma}
      
      \parag{Proof} As we know that \ $ \tilde{\mu}_{min}(G^b) =  \tilde{\mu}_{min}(G)$ , for \ $M \in \mathbb{N}$, the assumption \ $M > \mu -  \tilde{\mu}_{min}(G)$ \  implies that \ $(a - (\mu-M).b).G^b = b.G^b$ \ (see the remark before proposition \ref{sub min}). By  definition of the index of \ $G$ \ we have  \ $b^{\delta(G)}.G \subset G^b $. Combining both gives 
      $$ b^{M+\delta(G)+1}.G \subset b^M.(a -(\mu-M).b).G = (a-\mu.b).b^M.G \subset (a-\mu.b).G.$$
      Now let \ $N = M + \delta(G) +1$ ; a sufficient condition on the integer \ $N$ \  is now \ $N \geq \mu -\tilde{\mu}_{min}(G) + \delta(G) + 2.$ \ $\hfill \blacksquare$

     \begin{cor}\label{3}
     A sufficient condition for \ $N \geq N(F^*,-\lambda)$ \ in the situation of prop. \ref{4} in the regular case is that \ $N \geq or(E) + L(E) + rank(E) + 1 $.
     \end{cor}
     
     Remark that the inequality \ $L(E) + rank(E) \geq 1$ \ for any non zero regular \ $E$ \ implies that we have  \ $or(E) + L(E) + rank(E) + 1 \geq or(E)  +  2$.
     
     \parag{Proof} We apply the previous lemma with \ $ F^*= G$ \ and \ $\mu = -\lambda = -\tilde{\lambda}_{min}(E)$. The conclusion comes now from the following facts :
     \begin{enumerate}
     \item \ $-\widetilde{(-\lambda)}_{min}(F^*) = \tilde{\lambda}_{max}(F) \leq \lambda + L(E) $ \ this last inequality is proved in \ref{Induction largeur}.
          \item \ $ \delta(F^*) = \delta(F) \leq or(F) \leq or(E) $ \ proved in \ref{or dual} and \ref{ordre reg.} $\hfill \blacksquare$
          \end{enumerate}

     \subsection{The theorem.}
     
     \begin{thm}\label{finite det. thm}
     Let \ $E$ \ be a regular (a,b)-module. There exists an integer \ $N(E)$ \ such that for any (a,b)-module \ $E'$,  any integer \ $ N \geq N(E)$ \ and any \ $\A-$isomorphism
     \begin{equation*}
      \varphi :  E/b^N.E \to  E'/b^N.E'  \tag{1}
      \end{equation*}
     there exists an unique \ $\A-$isomorphism \ $ \Phi :  E \to E' $ \ inducing the given \ $\varphi$.\\
     Moreover the choice \ $N(E)= N_0(E) : =  or(E) + L(E) + rank(E) + 1 $ \ is possible.
     \end{thm}
     
     \parag{Remarks}
     \begin{enumerate}
     \item It is easy to see that for a rank 1 regular (a,b)-module the integer 2  is the best possible.
     \item In our final lemma \ref{final} we show that the integer given in the theorem is optimal for the rank $k$ (a,b)-module \ $J_k(\lambda)$, (defined in the lemma), for any \ $k \in \mathbb{N}^*$.
     \item For the rank 2 (a,b)-modules \ $E_{\lambda, \mu}$ \ the integer  given by the theorem  is \ $or(E) + L(E) + 2 + 1  =  3$ \ is again optimal, as it can be shown in the same maner that in our final lemma. 
      \item  For the rank 2 simple pole (a,b)-module \ $E_{\lambda}(0)$ \ the integer given by the theorem is  \ $ 3 = L(E) + rank(E) + 1 $ \ and  the best possible is \ $2$ : the action of \ $b^{-1}.a$ \ on \ $E/b.E$ \ which is determined by \ $E/b^2.E$ \ characterizes this rank 2 regular (a,b)-module in the classification given in proposition \ref{class}.  
         \item For the (a,b)-module \ $E$ \  associated to an holomorphic germ at the origine of \ $\mathbb{C}^{n+1}$ \ with an isolated singularity we have the uniform bounds \ $or(E) \leq n+1$ \ and \ $\ L(E) \leq n $ \ so the choice \ $N(E) = 2n+ \mu +2$ \ is always possible, where \ $\mu$ \ is the Milnor number (equal to the rank). \ $\hfill \square$
     \end{enumerate}
     
     \parag{Proof} We shall make an induction on the rank of \ $E$. So we shall assume that the result is proved for a rank \ $p - 1$ \ (a,b)-module and we shall consider a regular (a,b)-module \ $E$ \ of rank \ $p \geq 1$, an (a,b)-module \ $E'$, an integer \ $N \geq N_0(E)$ \ and an \ $\A-$isomorphism \ $\varphi$ \ as in \ $(1)$. From \ref{det. finie ordre reg.} we know that \ $E'$ \ is then regular and has order of regularity \ $or(E') = or(E)$.\\
         Choose now a complex number \ $\lambda$ \ which is minimal \ $modulo \ \mathbb{Z} $ \ such there exists an exact sequence of (a,b)-module ( so \ $\lambda = \tilde{\lambda}_{min}(E)$ \ with the terminology of \S 1.3)
     \begin{equation*}
      0 \to E_{\lambda} \overset{\alpha}{\to}  E \overset{\beta}{\to} F \to 0 . \tag{2}
      \end{equation*}
     This exists from theorem \ref{J-H}. The (a,b)-module \ $F$ \ has rank \ $p-1$ \ and from \ref{3} \ we have \ $N_0(E) \geq N(F^*, -\lambda)$. So we know from \ref{4} that the extension \ $(2)$ \ is determined by the extension
       \begin{equation*}
      0 \to E_{\lambda}/b^N.E_{\lambda} \overset{\alpha_N}{\to}  E/b^N.E \overset{\beta_N}{\to} F/b^N.F \to 0 . \tag{$2_N$}
      \end{equation*}
      Now, using the \ $\A-$isomorphism \ $\varphi$ \ we obain an injective \ $\A-$linear map
      $$ j_N :  E_{\lambda}/b^N.E_{\lambda} \hookrightarrow E'/b^N.E' .$$
      Using the proposition \ref{sub min} with the fact that \ $N \geq or(E') + 2$ \ there exists a unique normal inclusion \ $j : E_{\lambda} \hookrightarrow E' $ \ inducing \ $j_N$.\\
      Define \ $F' : = E'/j(E_{\lambda})$. Then \ $F'$ \ is a rank \ $p-1$ \ (a,b)-module and the exact sequence
      \begin{equation*}
      0 \to E_{\lambda} \overset{j}{\to} E' \to F' \to 0 \tag{2'}
      \end{equation*}
      induced the extension \ $(2_N)$. Using the induction hypothesis, because the inequalities \ $or(E) \geq or(F)$ \ from \ref{det. finie ordre reg.}  and \ $L(E) + rank(E) \geq L(F) + rank(F)$ \ from \ref{maj. lambda(max)(F)} implies \ $N_0(E) \geq N_0(F)$ , we have a unique isomorphism  \ $\Psi : F \to F' $ \ compatible with the one induced by \ $\varphi$ \ between \ $F/b^N.F$ \ and \ $F'/b^N.F'$. Using \ref{4},  \ref{3} and  the inequality \ $N_0(E) \geq N(F^*, -\lambda)$ \  we have an unique isomorphism of extensions
      $$\xymatrix{ 0 \ar[r] & E_{\lambda} \ar[d]^= \ar[r]^{\alpha} & E \ar[d]^{\Phi} \ar[r]^{\beta} & F \ar[d]^{\Psi} \ar[r] & 0 \\
      0 \ar[r] &  E_{\lambda} \ar[r]^j & E' \ar[r] & F' \ar[r] & 0 }$$
      concluding the proof. $\hfill \blacksquare$
      
      \bigskip
      
         \begin{lemma}\label{final} 
         Let \ $E : = J_k(\lambda) $ \ the rank \ $k$ \ (a,b)-module defined by the \ $\mathbb{C}[[b]]-$basis \ $e_1, \cdots, e_k$ \ and by the following relations
      $$ a.e_j = (\lambda + j - 1).b.e_j + e_{j+1} \quad \forall  j \in [1,k]  $$
      with the convention \ $e_{k+1} = 0 $. We have \ $\delta(E) = or(E) = k -1, L(E) = -k+1$. The integer \ $or(E) + L(E) + rank(E) + 1 = k+1$ \ is the best possible for the theorem.
      \end{lemma}
      
      \parag{Proof} It is easy to see that the saturation \ $E^{\sharp}$ \ is generated by \ $e_1, b^{-1}.e_2, \cdots, b^{-k+1}.e_k$. This gives the equality \ $\delta(E) = or(E) = k -1$.\\
      Assume that we have an inclusion \ $E_{\mu} \hookrightarrow E$ \ such that \ $e_{\mu} \not\in b.E$. Then there exists \ $(\alpha_1, \cdots, \alpha_k) \in \mathbb{C}^k\setminus\{0\} $ \ such that
      $$ a.(\sum_{h =1}^k \ \alpha_h.e_h) = \mu.b(\sum_{h=1}^k \alpha_h.e_h) + b^2.E .$$
      Then we obtain
      $$ \sum_{h=1}^k \ \alpha_h.\big((\lambda + j-1).b.e_h + e_{h+1}\big) =  \sum_{h=1}^k \ \alpha_h.\mu.b.e_h + b^2.E $$
      and so \ $\alpha_1=  \cdots = \alpha_{k-1} = 0 $ \ and we conclude that \ $\mu = \lambda + k - 1$.\\
      An easy computation shows that \ $J_k(\lambda)^* = J_k(-\lambda-2k +2)$ \ and so we have \ $\lambda_{max} = \lambda$. So \ $L(E) = -k + 1$.\\
      Now we shall prove that the integer \ $k + 1$ \ is the best possible in the theorem \ref{finite det. thm} for \ $E = J_k(\lambda)$ \ by giving a regular (a,b)-module  \ $F$ \ such that \ $F/b^k.F \simeq E/b^k.E$ \ and not isomorphic to \ $E$.\\
      Let consider the rank \ $k$ \ (a,b)-module \ $F$ \ defined by  \ $ \sum_{j=1}^k \ \mathbb{C}[[b]].e_j$ \ with the following relations
      \begin{align*}
      & a.e_j =  (\lambda + j -1).b.e_j + e_{j+1} \quad \forall j \in [1,k] \\
      & a.e_k = (\lambda + k -1).b.e_k + \sum_{h=1}^{k-1} \ \alpha_h.b^{k-h+1}.e_h 
      \end{align*}
      Then define, for \ $\beta_1, \cdots, \beta_{k-1} \in \mathbb{C}$,
      $$\varepsilon : = e_k + \sum_{j=1}^{k-1} \ \beta_j.b^{k-j}.e_j .$$
      We have
   \begin{align*}
     a.\varepsilon : = &(\lambda + k -1).b.e_k  + \sum_{h=1}^{k-1} \ \alpha_h.b^{k-h+1}.e_h  \ +\\
    \qquad  \qquad   & \sum_{j=1}^{k-1} \beta_j.\big[b^{k-j}.\big((\lambda+j-1).b.e_j + e_{j+1}) + (k-j).b^{k-j+1}.e_j \big] \\
    a.\varepsilon : = &(\lambda + k -1).b.e_k  + \sum_{h=1}^{k-1} \big(\alpha_h + \beta_h.(\lambda + k-1) + \beta_{h-1}\big).b^{k-h+1}.e_h 
   \end{align*}
   
   Let now choose \ $\beta_1, \cdots, \beta_{k-1}$ \ such that we have 
   $$ \alpha_h + \beta_h.(\lambda + k - 1) + \beta_{h-1} = (\lambda + k -1 + \beta_{k-1}).\beta_h \quad \forall h \in [1,k-1]$$
   with the convention \ $\beta_0 = 0 $. We obtain the system of equations
   $$ \alpha_h + \beta_{h-1} = \beta_{k-1}.\beta_h \quad \forall h \in [1,k-1].$$
   This implies, assuming \ $\beta_{k-1} \not= 0$, that \ $\beta_{k-1}$ \ is solution of the equation
   $$ x^k = \alpha_{k-1}.x^{k-2} + \cdots + \alpha_2.x + \alpha_1.$$
   Now choose \ $\alpha_2 = \cdots = \alpha_{k-1} = 0 $ \ and \ $\alpha_1 : = \rho^k $ \ with \ $\rho \in ]0,1[$. Then  choose \ $\beta_j = \rho^{k-j} \quad \forall j \in [1,k-1]$. It is clear that the corresponding \ $F_{\rho}$ \ satisfies \ $F/b^k.F \simeq E/b^k.E$ \ as \ $a.e_k  = e_k + \rho.b^k.e_1$ \ in \ $F_{\rho}$. But the relation \ $a.\varepsilon = (\lambda + k -1 + \rho^k).b.\varepsilon $ \ with \ $\varepsilon \not= 0$ \  shows that \ $F_{\rho}$ \ cannot be isomorphic to \ $J_k(\lambda)$. $\hfill \blacksquare$   
   
   \newpage

      \section{Bibliography}
      
      \bigskip
      
      \begin{enumerate}
      
      \item{[B.93]} Barlet, Daniel \textit{Theory of (a,b)-modules I} \ in Complex Analysis and Geometry, Plenum Press New York (1993), p.1-43.
      
      \item{[B.95]} Barlet, Daniel \textit{Theorie des (a,b)-modules II. Extensions} in Complex Analysis and Geometry, Pitman Research Notes in Mathematics Series  366 Longman (1997), p. 19-59.
      
      \item{[B.07]} Barlet, Daniel {\it Sur certaines singularit\'es non isol\'ees d'hypersurfaces II}, to appear in  the Journal of Algebraic Geometry. Preprint Inst. E. Cartan (Nancy) (2006) n. 34 59 pages.
      
      \item{[D.70]} Deligne, Pierre \textit{\'Equations diff\'erentielles \`a points singuliers reguliers} Lect. Notes in Maths, vol. 163, Springer-Verlag (1970).
      \item{[Gr.66]} Grothendieck, A. {\it On the de Rham cohomology of algebraic varieties}, Publ.Math. IHES 29 (1966), p. 93-101.
      \item{[M.91]} Malgrange, Bernard {\it Equations diff\'erentielles \`a coefficients polynomiaux}, Progress in Mathematics  vol. 96 (1991) Birkh{\"a}user Boston.
      
      \item{[Sa.91]} Saito, Morihiko {\it Period mapping via Brieskorn modules}, Bull. Soc. Math. France vol. 119 n.2 (1991), p. 141-171.
      
      \item{[St.76]} Steenbrink, Joseph {\it Mixed Hodge structures on th vanishing cohomology}, Proc. Nordic Summer School on Real and Complex Singularities, Oslo (1976) Sijthoff and Noordhoff 1977.
      
      \item{[Va.81]}Varchenko, A. N. {\it An asymptotic  mixed Hodge structure in vanishing cohomologies}, Izv. Acad. Sci. SSSR, ser. mat. 45 (1981), p. 540-591.
      
      \end{enumerate}

 \end{document}